\documentclass[12pt]{article}

\usepackage{amsmath}
\usepackage{amssymb}
\usepackage{amsthm}
\usepackage{latexsym}
\usepackage{color}
\usepackage{graphicx}
\usepackage{color}

\DeclareSymbolFont{calletters}{OMS}{cmsy}{m}{n}
\DeclareSymbolFontAlphabet{\mathcal}{calletters}

\def\be{\begin{eqnarray}}
\def\ee{\end{eqnarray}}

\def\b*{\begin{eqnarray*}}
\def\e*{\end{eqnarray*}}

\newtheorem{Theorem}{Theorem}[section]
\newtheorem{Definition}[Theorem]{Definition}
\newtheorem{Proposition}[Theorem]{Proposition}

\newtheorem{Assumption}[Theorem]{Assumption}

\newtheorem{Lemma}[Theorem]{Lemma}

\newtheorem{Remark}[Theorem]{Remark}
\newtheorem{Example}[Theorem]{Example}

\makeatletter \@addtoreset{equation}{section}

\usepackage{color}

\newcommand{\cblue}{\color{blue}}

\def \E{\mathbb{E}}
\def \F{\mathbb{F}}
\def \L{\mathbb{L}}
\def \P{\mathbb{P}}

\def \R{\mathbb{R}}
\def \S{\mathbb{S}}

\def \N{\mathbb{N}}

\def\Ac{{\cal A}}

\def\Ec{{\cal E}}
\def\Fc{{\cal F}}
\def\Gc{{\cal G}}

\def\Kc{{\cal K}}
\def\Lc{{\cal L}}
\def\Mc{{\cal M}}

\def\Oc{{\cal O}}

\def\Tc{{\cal T}}
\def\Uc{{\cal U}}

\def\Wc{{\cal W}}

\def\Yc{{\cal Y}}
\def\Zc{{\cal Z}}

\def\Cb{\overline{C}}
\def\Cu{\underline{C}}
\def\Ob{\overline{O}}
\def\Ocb{\overline{\Oc}}
\def\Bh{\widehat{B}}
\def\Lcb{\overline{\Lc}}

\def\Nt{\tilde{N}}
\def\Bt{\tilde{B}}
\def\Kct{\tilde{\Kc}}
\def\Omt{\tilde{\Om}}
\def\Fct{\tilde{\Fc}}
\def\Pt{\tilde{\P}}
\def\pt{\tilde{p}}
\def\Yt{\tilde{Y}}
\def\Zt{\tilde{Z}}

\def\Acu{\overline{\Ac}}
\def\Acb{\underline{\Ac}}
\def\Ucu{\overline{\Uc}}
\def\Ucb{\underline{\Uc}}

\def\taub{\bar \tau}

\def\thetau{\bar{\theta}}
\def\psiu{\overline{\psi}}
\def\psib{\underline{\psi}}

\def \Om{\Omega}
\def \om{\omega}

\def \omh{\hat{\om}}
\def \tauh{\hat{\tau}}
\def \omt{\tilde{\om}}
\def \eps{\varepsilon}
\def \xb{\mathbf{x}}
\def \xbh{\hat{\xb}}
\def \0{\mathbf{0}}
\def \H{\textsc{h}}

\def \Fcb{\overline{{\cal F}}}
\def \Fbb{\overline{\mathbb{F}}}

\newcommand{\rmi}{{\rm (i)$\>\>$}}

\newcommand{\rmii}{{\rm (ii)$\>\>$}}
\newcommand{\rmiii}{{\rm (iii)$\>\,    \,$}}

\newcommand{\rma}{{\rm a)$\>\>$}}
\newcommand{\rmb}{{\rm b)$\>\>$}}
\newcommand{\rmc}{{\rm c)$\>\>$}}

\def\x{\times}
\def\ox{\otimes}
\def\1{{\bf 1}}
\def \proof{{\noindent \bf Proof. }}

\title{A numerical algorithm for a class of BSDEs via branching process \footnote{The authors are grateful to Emmanuel Gobet, Jia Zhuo and two anonymous referees for very helpful comments and suggestions.}}

\author{Pierre Henry-Labord\`ere\thanks{Soci\'et\'e G\'en\'erale, Global Market Quantitative Research,
pierre.henry-labordere@sgcib.com}
        \and Xiaolu Tan\thanks{CEREMADE, University of Paris-Dauphine, xiaolu.tan@ceremade.dauphine.fr}
        \and Nizar Touzi\thanks{Ecole Polytechnique Paris, Centre de Math\'ematiques Appliqu\'ees,
        nizar.touzi@polytechnique.edu. Research supported by the Chair {\it Financial Risks} of the {\it Risk Foundation} sponsored by Soci\'et\'e G\'en\'erale, and the Chair {\it Finance and Sustainable
Development} sponsored by EDF and Calyon.} }

\date{\today}

\begin{document}

\maketitle

\abstract{
        We give a study to the algorithm for semi-linear parabolic PDEs in Henry-Labord\`ere
        \cite{Henry-Labordere_branching} and then generalize it to the non-Markovian case for a class of Backward SDEs (BSDEs).
        By simulating the branching process, the algorithm does not need any backward regression.
        To prove that the numerical algorithm converges to the solution of BSDEs, we use the notion
        of viscosity solution of path dependent PDEs introduced by Ekren, Keller, Touzi and Zhang
        \cite{EkrenKellerTouziZhang} and extended in Ekren, Touzi and Zhang
        \cite{EkrenTouziZhang1, EkrenTouziZhang2}.

\vspace{1mm}

 }

\section{Introduction}

        Initially proposed by Pardoux and Peng \cite{PardouxPeng}, the theory of Backward Stochastic Differential
        Equation (BSDE) has been largely developed and has many applications in control theory, finance etc.
        In particular, BSDEs can be seen as providing a nonlinear Feynman-Kac
        formula for semi-linear parabolic PDEs in the Markovian case, i.e. the solution
        of a Markovian type BSDE can be given as the
        viscosity solution of a semi-linear PDE. We also remark that this connection has been  extended recently to
        the non-Markovian case by Ekren, Keller, Touzi and Zhang \cite{EkrenKellerTouziZhang} with the notion of viscosity
        solution of  path dependent PDEs (PPDEs).

        Numerical methods for BSDE have also been largely investigated since then.
        The classical numerical schemes for BSDEs are usually given as a backward iteration, where every step
        consists in estimating the conditional expectations, see e.g. Bouchard and Touzi \cite{BouchardTouzi},
        Zhang \cite{Zhang}.
        Generally, we use the regression method to compute the conditional expectations, which is quite costly
        in practice and suffers from the curse of dimensionality.

        Recently, a new numerical algorithm has been proposed by Henry-Labord\`ere \cite{Henry-Labordere_branching} for a class of semi-linear PDEs, using an extension of branching process.
        First, it is a classical result that the branching diffusion process gives a probabilistic representation of
        the so-called KPP (Kolmogorov-Petrovskii-Piskunov) semi-linear
        PDE (see e.g. Watanabe \cite{Watanabe}, McKean
        \cite{McKean_1975}):
        \be \label{eq:PDEsemiLinear}
            \partial_t u(t,x) ~+~ \frac{1}{2} D^2 u(t,x) ~+~ \beta \Big( \sum_{k=0}^{n_0} a_k u^k(t,x) ~-~ u(t,x) \Big)
            &=& 0
        \ee
        with a terminal condition $u(T,x) = \psi(x)$, where $D^2$ is the Laplacian on $\R^d$, and $(a_k)_{0 \le k \le n_0}$ is a
        probability mass sequence, i.e. $a_k \ge 0$ and $\sum_{k=0}^{n_0} a_k = 1$. The above
        semi-linear PDE \eqref{eq:PDEsemiLinear} characterizes a branching Brownian motion,
        where every particle in the system dies in
        an exponential time of parameter $\beta$, and creates $k$ i.i.d. descendants with probability $a_k$.
        More precisely, let $N_T$ denote the number of particles alive at time $T$, and $(Z_T^i)_{i=1, \cdots, N_T}$
        denote the position of each particle, then up to integrability, the function
        \b*
            v(t,x) &:=& \E_{t,x} ~ \big[ \Pi_{i=1}^{N_T} ~ \psi(Z_T^i) \big]
        \e*
        solves the above equation \eqref{eq:PDEsemiLinear}, where the subscript ${t,x}$ means that the system is
        started at time $t$ with one particle at position $x$.
        This connection has then also been extended for a larger class of nonlinearity,
        typically $u^\alpha, \alpha \in [0,2]$, with the superdiffusion,
        for which we refer to Dynkin \cite{Dynkin} and Etheridge \cite{Etheridge}.
		Moreover, this representation allows to solve numerically the PDE
        \eqref{eq:PDEsemiLinear} by simulating the corresponding branching
        process.

        When the coefficients $a_k$ are arbitrary in equation \eqref{eq:PDEsemiLinear} and
        the Laplacian $D^2$ is replaced by an It\^o operator $\Lc_0$ of the form
        \b*
            \Lc_0 u(t,x) &:=& \mu(t,x) \cdot D u(t,x) ~+~ \frac{1}{2} \sigma \sigma^T(t,x) : D^2 u(t,x),
        \e*
        Henry-Labord\`ere's \cite{Henry-Labordere_branching} proposed to simulate a branching diffusion process
        with a probability mass sequence $(p_k)_{k=0,\cdots,M}$, and by counting the weight
$\frac{a_k}{p_k}$, he obtained a
        so-called ``marked'' branching diffusion method.
        A sufficient condition for the convergence of the algorithm is provided
        in \cite{Henry-Labordere_branching}.
        In particular, the algorithm does not need to use the regression method,
        which is one of the main advantages comparing to the BSDE method.

        For PDEs of the form \eqref{eq:PDEsemiLinear}, Rasulov, Raimov and Mascagni
        \cite{RasulovRaimovMascagni} introduced
        also a Monte-Carlo method using branching processes. Their method depends essentially on
        the representation of its solution by the fundamental solution of the heat equation.

        The main objective of this paper is to give a more rigorous study to the algorithm in
        \cite{Henry-Labordere_branching} and also to  generalize it to the non-Markovian case for a class of decoupled Forward
        Backward SDEs (FBSDEs) whose generators can be represented as the sum of a power series,
        which can be formally approximated by polynomials.
        Although the polynomial generators are only locally Lipschitz, the solutions may be
        uniformly bounded under appropriate conditions, and hence they can be considered as standard
        decoupled FBSDEs with Lipschitz generators.

        Our numerical solution is based on a branching process, which is constructed by countably many independent Brownian motions and exponential random variables.
        To bring back the numerical solution to the BSDE context of one Brownian motion with the Brownian filtration, we use the notion of viscosity solution of path dependent PDEs introduced by Ekren, Keller, Touzi and Zhang \cite{EkrenKellerTouziZhang} and next extended in Ekren, Touzi and Zhang \cite{EkrenTouziZhang1, EkrenTouziZhang2}.
        Namely, we shall prove that the numerical solution obtained by the branching diffusion  is a viscosity solution to a corresponding semilinear PPDE, which admits also a representation by a decoupled FBSDE as illustrated in \cite{EkrenKellerTouziZhang}.
        Then the numerical solution is the unique solution of the corresponding FBSDE by the uniqueness of the solution to PPDEs.

        \vspace{2mm}

        The rest of the paper is organized as follows.
        In Section \ref{sec:BranchFBSDE}, we consider a class of decoupled FBSDEs whose generators can be represented as a convergent power series.
        We then introduce a branching diffusion process, which gives a representation of the solution of the FBSDE with polynomial generator.
        In particular, such a representation induces a numerical algorithm for the class of FBSDEs using branching process.
        Then in Section \ref{sec:proof_property_v}, we complete the proof of the regularity property of the value function represented by branching process.
        Next, we complete the proof of the main representation theorem in Section \ref{sec:proof_mainthm}.
        For this purpose, we introduce in Section \ref{subsec:PPDE} a notion of viscosity solution to a class of semilinear PPDE, where there is no non-linearity on the derivatives of the solution function, following Ekren, Touzi and Zhang \cite{EkrenTouziZhang1, EkrenTouziZhang2}.
        The uniqueness of solution to our PPDE and its representation by FBSDE are proved by the same arguments as in \cite{EkrenTouziZhang1, EkrenTouziZhang2}, which are hence provided in Appendix.
        Finally, we illustrate the efficiency of our algorithm by some numerical examples in Section \ref{sec:num_example}.

\section{A numerical algorithm for a class of BSDEs}
\label{sec:BranchFBSDE}

        In this section, we shall consider a class of decoupled FBSDEs whose generators can be represented as a convergent power series, which can be approximated by polynomials.
        Then for FBSDEs with polynomial generators, we provide a representation of their solutions by branching diffusion processes.
        In particular, the representation induces a natural numerical algorithm for the class of FBSDEs by simulating the branching diffusion process.

\subsection{A class of decoupled FBSDEs}

Let $\Om^0 := \big\{ \om \in C([0,T], \R^d) ~: \om_0 = \0 \big \}$
be the canonical space of continuous paths with initial value $\0$,
$\F^0$ the canonical filtration and $\Lambda^0 := [0,T] \x \Om^0$.
For every $(t,\om) \in \Lambda^0$, denote $\| \om \|_t := \sup_{0
\le s \le t} |\om(s)|$.

        Then the canonical process $B(\omega)=\{B_t(\om) := \om_t, 0\le t\le T\}$ for all $\om \in \Om^0$, defines a Brownian motion under the Wiener measure $\P_0$.

        Let $\mu : \Lambda^0 \to \R^d$ and $\sigma : \Lambda^0 \to \S^d$ be $\F^0-$progressively measurable processes.
        Suppose further that for every $0 \le t \le t' \le T$ and $\om, \om' \in \Om^0$,
        \be \label{eq:condi_mu_sigma}
            | \mu(t,\om) - \mu(t', \om')| + | \sigma(t, \om) - \sigma(t', \om')|
            &\le&
            C \big(\sqrt{|t - t'|} + \| \om_{t \land \cdot} - \om'_{t'\land \cdot} \|_T  \big)
        \ee
        for some constant $C >0$, and $\sigma \sigma^T(t,\om) \ge c_0 I_d$ for some constant $c_0 > 0$. We denote, for every $(t,\xb) \in \Lambda^0$, by $^{t,\xb}X$ the solution of the following SDE under $\P_0$:
        \be \label{eq:SDE}
            X_s = \xb_s,~ \forall s \le t
            ~&\mbox{and}&
            X_s = \xb_t + \int_t^s \mu(r,X_{\cdot}) dr + \int_t^s \sigma(r,X_{\cdot}) dB_r,~ \forall s > t.
        \ee
        For later uses, we provide an estimate on the SDE \eqref{eq:SDE}.
        \begin{Lemma} \label{lemm:SDE_Lip}
            There is a constant $C$ such that for every $t \in [0,T]$ and $(t_1, \xb_1), (t_2, \xb_2) \in [t,T] \x \Om^0$,
            \b*
                \E^{\P_0} \Big[ \sup_{t \le s \le T}
                    \Big| ~^{t, \xb_1}X_{s \land t_1} - ~^{t, \xb_2}X_{s \land t_2} \Big|^2
                    \Big]
                \le
                C \big( 1 +\|\xb_1\|_{t}^2  +
                \| \xb_2 \|_{t}^2  \big) \big( |t_1 - t_2| + \| \xb_1 - \xb_2\|_t^2 \big).
            \e*
        \end{Lemma}
        \proof Suppose, without loss of generality, that $t_1 \le t_2$, we notice that
        \b*
            \E^{\P_0} \Big[ \sup_{t \le s \le T}
                    \Big| ~^{t, \xb_1}X_{s \land t_1} - ~^{t, \xb_2}X_{s \land t_2} \Big|^2
                    \Big]
            &\le&
            \E^{\P_0} \Big[ \sup_{t \le s \le t_1}
                    \Big| ~^{t, \xb_1}X_s - ~^{t, \xb_2}X_s \Big|^2
                    \Big] \\
            &&+~
            \E^{\P_0} \Big[ \sup_{t_1 \le s \le t_2}
                    \Big| ~^{t, \xb_2}X_{t_1} - ~^{t, \xb_2}X_{s} \Big|^2
                    \Big].
        \e*
        Then the estimate in Lemma \ref{lemm:SDE_Lip} is a standard result for SDEs, by using It\^o's formula and Gronwall's Lemma.
        One can find the arguments in Lemma 2 and Theorem 37 in Chapter V of Protter \cite{Protter} for an almost the same result. \qed

        \vspace{2mm}

        Suppose that $\psi: \Om^0 \to \R$ is a non-zero, bounded Lipschitz continuous function, and $F: (t,\xb,y) \in \Lambda^0 \x \R \to \R$ is a function Lipschitz in $y$ such that for every $y$, $F(\cdot, y)$ defined on $\Lambda^0$ is $\F^0-$progressive.
        We consider the following BSDE:
        \be \label{eq:FBSDE_poly}
            Y_t
            &=&
            \psi(~^{0,\0}X_{\cdot})
            +
            \int_t^T F(s,~^{0,\0}X_{\cdot}, Y_s)  ds
            - \int_t^T Z_s dB_s, ~~ \P_0-a.s.,
        \ee
    where the generator $F$ has the following power series representation in $y$, locally in $(t,\xb)$:
        \be \label{eq:def_f}
            F(t, \xb, y)
            &:=&
            \beta \Big( \sum_{k=0}^{\infty} a_k(t,\xb) y^k
                        ~-~ y
                  \Big),
            ~~ (t,\xb)\in\Lambda^0,
        \ee
for some constant $\beta > 0$, and some sequence $(a_k)_{k \ge 0}$
of bounded scalar $\F^0-$progressive functions defined $\Lambda^0$.
We also assume that every $a_k$ is uniformly $1/2-$H\"older-continuous in $t$ and
Lipschitz-continuous in $\om$.

Denoting by $|.|_0$ the $\L^\infty(\Lambda^0)$-norm, we
now formulate conditions on the power series
 \be\label{eq:def_ell}
 \ell_0(s)
 :=
 \sum_{k\ge 0}|a_k|_0\; s^k
 &\mbox{and}&
 \ell(s):=\beta\big[|\psi|_0^{-1}\ell_0(s|\psi|_0)-s\big],
 ~~s\ge 0,
 \ee
so as to ensure the existence and uniqueness of the solution to
BSDE \eqref{eq:FBSDE_poly} (see also Remark \ref{rem:conditionl} for an intuitive interpretation of the condition).

\begin{Assumption} \label{assum:function_p}
        \rmi The power series $\ell_0$ has a radius of convergence $0 < R \le \infty$, i.e. $\ell_0(s)<\infty$ for $|s|<R$ and
        $\ell_0(s)=\infty$ for $|s|>R$. Moreover, the function $\ell$ satisfies either one of the following
        conditions:
        \\
        $~\hspace{5mm}(\ell 1)$\quad $\ell(1) \le 0$,
        \\
        $~\hspace{5mm}(\ell 2)$\quad or, $\ell(1) > 0$ and for some $\hat s > 1$, $\ell(s) > 0, \forall s \in [1, \hat s)$ and $\ell(\hat s)=0$.
        \\
        $~\hspace{5mm}(\ell 3)$\quad or, $\ell(s)>0, \forall s \in [1, \infty)$ and $\int_1^{\bar s} \frac{1}{\ell(s)} ds=T,$ for some
        constant $s \in (1, \frac{R}{|\psi|_0})$.
        \\
        \rmii The terminal function satisfies $|\psi|_0 < R$.
\end{Assumption}

\begin{Proposition} \label{prop:Lip_BSDE}
Let Assumption \ref{assum:function_p} hold true, then the BSDE
\eqref{eq:FBSDE_poly} has a unique solution $(Y,Z)$ such that
$\sup_{0 \le t \le T} |Y_t| \le R_0$, $\P_0-$almost surely for some constant $R_0 < R$.
\end{Proposition}

\begin{Remark}
            When $\psi \equiv 0$, the function $\ell$ in \eqref{eq:def_ell} is not well defined. In order to provide a sufficient condition for the power series representation, we can consider the BSDE \eqref{eq:FBSDE_poly} with terminal condition $Y_T = \eps$.
            Define the corresponding function $q_{\eps}(s) := \beta \big[\eps^{-1}\ell_0(\eps s)-s \big] $.
            Suppose that for some $\eps > 0$, Assumption \ref{assum:function_p}
            holds true with the corresponding function $q_{\eps}$, then by comparison result of standard BSDEs with global Lipschitz generator, the BSDE \eqref{eq:FBSDE_poly} admits still a unique solution $(Y,Z)$ such that $Y$ is uniformly bounded (notice that when $Y$ is uniformly bounded, the generator $F$ is Lipschitz in $y$).
        \end{Remark}

        In preparation of the proof, let us consider first the ordinary differential equation (ODE) of $\rho(t)$ on interval $[0,T]$:
        \be \label{eq:ODE}
             \rho' ~=~ \ell(\rho), &\mbox{with initial condition}& \rho(0) ~=~ 1.
        \ee

        \begin{Lemma} \label{lemm:solution_ODE}
            Let $|\psi|_0 < R$, then ODE \eqref{eq:ODE} admits a unique bounded solution on the interval
            $[0,T]$ if and only if Assumption \ref{assum:function_p} \rmi holds true.
          	Moreover, in this case, we have $ 0 \le \rho(t) \le \frac{R_0}{|\psi|_0},~ \forall t \in [0,T]$ for some constant $R_0 < R$.
        \end{Lemma}
        \proof First, since the function $\ell$ is Lipschitz on $[0, L]$ for every $L < \frac{R}{|\psi|_0} $, then
        it follows by Picard-Lindel\"of theorem (see e.g. Chapter 2 of Teschl \cite{Teschl}) that there is $T_{\max} > 0$ such that ODE \eqref{eq:ODE} admits a unique solution $\rho$ on $[0, T_{\max})$ and that $\lim_{t \to T_{\max}} |\rho(t)| = \frac{R}{|\psi|_0} > 1$.
        Further, we observe that $\ell(0) \ge 0$, which implies that $\rho(t) \ge 0$ on $[0, T_{\max})$.
        Then it is enough to prove that $T_{\max} > T$.

        Let us now discuss three cases of Assumption \ref{assum:function_p}.
        \rmi Suppose that $(\ell 1)$ holds true, i.e $\ell(1) \le 0$. It follows then $\rho(t) \in [0,1]$ for every $t \in [0, T_{\max})$ and hence $T_{\max} = \infty > T$.
        \rmii Suppose now $\ell(1) > 0$ and for some $\hat s > 1$, $\ell(s) > 0, \forall s \in [1, \hat s)$ and $\ell (\hat s) =0$.
        It is clear that in this case, $t \mapsto \rho(t)$ is increasing and $\rho(t)$ converges to $\hat s$ as $t \to \infty$, and hence $T_{\max} = \infty > T$.
        \rmiii Otherwise, suppose that $(\ell 3)$ holds true, it follows then by \eqref{eq:ODE} that
        \b*
            T_{\max}
            =
            \int_0^{T_{\max}} dt
            =
            \int_0^{T_{\max}} \frac{1}{\ell(\rho(t))} d \rho(t)
            =
            \int_1^{R/|\psi|_0} \frac{1}{\ell(s)} ds,
        \e*
        since $\rho(0) = 1$ and $\rho(T_{\max}) = R/|\psi|_0$.
        We hence deduce that $T < T_{\max}$ by Assumption \ref{assum:function_p} \rmi $(\ell 3)$ and the positivity of the function $\ell$ on $[1, \infty)$.
        \qed

        \begin{Remark} \label{rem:eqv_ODE}
            The ODE \eqref{eq:ODE} can be rewritten as
            \b*
            \rho(t)
            &=&
            \rho(0)
            +\int_0^t \ell(\rho(s)) ds.
            \e*
            Let $\varphi(t) := \rho(t) |\psi|_0$, then under Assumption \ref{assum:function_p} we have
            \be \label{eq:ODE_v}
                    e^{\beta t} \varphi(t) &=& \varphi(0) ~+~ \int_0^t e^{\beta s} \beta \Big( \sum_{k=0}^{\infty} |a_k|_0 \varphi^k(s) \Big) ds.
            \ee
            In other words, the existence and uniqueness of solution to \eqref{eq:ODE} is equivalent to that of \eqref{eq:ODE_v}.
        \end{Remark}

        \begin{Remark} \label{rem:ODE_stability}
                Suppose that $a_k \equiv 0$ for every $k > n_0$ with some $n_0 \in \N$, then clearly
                $\ell (s) := \beta \Big(\sum_{k=0}^{n_0} |a_k|_0 |\psi|_0^{k-1} s^k - s \Big)$ and the convergence radius $R = \infty$.
                Denote $\ell_{\eps}(s) := \beta \big( \sum_{k=0}^{n_0} |(1+\eps)a_k|_0 |(1+\eps)\psi|_0^{k-1} s^k -s \big)$.
                Let Assumption \ref{assum:function_p} hold true for $\ell$, then for $\eps > 0$ small enough, 
                $\ell_{\eps}$ also satisfies one of the conditions $(\ell 1 - \ell 3)$ in Assumption \ref{assum:function_p}.
                It follows that the ODE: $\rho'(t) = \ell_{\eps}(\rho)$ with initial condition $\rho(0) = 1$ admits a unique solution on $[0,T]$ under Assumption \ref{assum:function_p}.
        \end{Remark}

    With the above existence and uniqueness result of ODE \eqref{eq:ODE}, we get the existence and uniqueness of the BSDE \eqref{eq:FBSDE_poly} in Proposition \ref{prop:Lip_BSDE}.

    \vspace{2mm}

        \noindent {\bf Proof of Proposition \ref{prop:Lip_BSDE}}. By Lemma \ref{lemm:solution_ODE},
        the solution $\rho$ of ODE \eqref{eq:ODE} is uniformly bounded by
        $\frac{C}{|\psi|_0}$ with some constant $C = R_0 < R$, where
        $R$ is the convergence radius of the power series $\sum_{k=0}^{\infty} |a_k|_0 x^k$.
        Denote $y_C := -C \vee (y \land C)$ for every $y \in \R$,
        \b*
        F_C(s, \xb, y)
        &:=&
        F(s, \xb, y_C)
        \e*
and
        \b*
        \overline{f}_C(s, \xb, y) ~:=~ \beta \Big( \sum_{k=0}^{\infty} |a_k|_0 |y_C|^k - y_C \Big),
            &&
            \underline{f}_C(s, \xb, y) ~:= - \beta \Big( \sum_{k=0}^{\infty} |a_k|_0 |y_C|^k + |y_C| \Big).
        \e*
Then $F_C$, $\overline{f}_C$ and $\underline{f}_C$ are all globally
Lipschitz in $y$, and $\underline{f}_C \le F_C \le \overline{f}_C$.
Moreover, if we replace the generator $F$ by $\overline{f}_C$(resp.
$\underline{f}_C$), and the terminal condition $\psi$ by $|\psi|_0$
(resp. $-|\psi|_0$) in BSDE \eqref{eq:FBSDE_poly}, the solution is
given by $\overline{Z} := 0$ (resp. $\underline{Z} := 0$) and
        \b*
            \overline{Y}_t ~:=~ \rho(T-t) |\psi|_0 &&\mbox{(resp.}~ \underline{Y}_t ~:=~ - \rho(T-t) |\psi|_0).
        \e*
        Therefore, by comparison principle, it follows that the solution $(Y_C, Z_C)$ of BSDE \eqref{eq:FBSDE_poly}
        with generator $f_C$ satisfies $ \underline{Y} \le Y_C \le  \overline{Y} $, and hence $|Y_C| \le C$.
        Further, since $F(t,\xb,y) = F_C(t,\xb,y)$ for all $|y| \le C$, it follows that $(Y_C, Z_C)$ is the required solution of BSDE \eqref{eq:FBSDE_poly}.  \qed

        \begin{Remark} \label{rem:conditionl}
            When $(a_k)_{k \ge 0}$ and $\psi$ are all positive constant functions, then the BSDE \eqref{eq:FBSDE_poly} degenerates to an ODE of the form \eqref{eq:ODE}.
            That is also the main reason for which we suppose Assumption \ref{assum:function_p} to guarantee the existence and uniqueness of the BSDE \eqref{eq:FBSDE_poly}.
        \end{Remark}

\subsection{A branching diffusion process}
\label{subsec:branch_diffusion}

        Let $\beta > 0$, $n_0 \ge 0$ and $p = (p_k)_{0 \le k \le {n_0}}$ be such that $\sum_{k\le n_0}p_k=1$ and $p_k \ge 0,~ k = 0, \cdots, {n_0}$. We now construct a branching diffusion process as follows: a particle starts at time $t$, from position $x$, performs a diffusion process given by \eqref{eq:SDE}, dies after a mean $\beta$ exponential time and produces $k$ i.i.d. descendants with probability $p_k$. Then the descendants go on to perform diffusion process defined by \eqref{eq:SDE} driven by independent Brownian motions. Every descendant dies and reproduces i.i.d. descendants independently after independent exponential times, etc. In the following, we shall give a mathematical construction of this branching diffusion process in three steps.

		In preparation, let $(\Om, \Fc, \P)$ be an abstract probability space containing a sequence of independent $d$-dimensional standard Brownian motions $(W^k)_{k \ge 1}$, a sequence of i.i.d. random variables  $(T^{i,j})_{i,j \ge 0}$ as well as i.i.d. r.v. $(I_n)_{n \ge 1}$, where $T^{0,0}$ is of exponential distribution $\Ec(\beta)$ with mean $\beta> 0$ and $I_1$ is of multi-nomial distribution $\Mc(p)$, i.e. $\P(I_1 = k) = p_k$, $\forall k = 0, 1 ,\cdots, {n_0}$. Moreover, the sequences $(W^k)_{k \ge 1}$,  $(T^{i,j})_{i,j \ge 0}$ and  $(I_n)_{n \ge 1}$ are mutually independent.

\paragraph{A birth-death process}
        We shall construct a continuous-time birth-death process associated with the coefficient $\beta > 0$ and the probability density sequence $(p_k)_{0 \le k \le {n_0}}$.

        The branching process starts with a particle at time $0$, $N_t$ denotes the number of the particles in the system, every particle runs an independent exponential time and then branches into $k$ i.i.d. particles with probability $p_k$.
        We denote by $T_n$ the $n-$th branching time of the whole system, at which one of the existing particles branches into $I_n$ particles.
        Between $T_n$ and $T_{n+1}$, every particle is indexed by $(k_1, \cdots, k_n) \in \{1, \cdots, {n_0}\}^n$, which means that its parent particle is indexed by $(k_1, \cdots, k_{n-1})$ between $T_{n-1}$ and $T_n$.
        We also have a bijection $c$ between $\N$ and $\cup_{n \ge 1} 2^{\{1, \cdots, {n_0} \}^n}$ defined by
        \be \label{eq:def_c}
        		c((k_1, \cdots, k_n)) ~:=~ \sum_{i=1}^n k_i ({n_0}+1)^i.
        \ee
        Denote by $\Kc_t$ the collection of the indexes of all existing particles in the system at time $t$. Then the initial setting of the system is given by
        \b*
                N_0 = 1,~~&T_0 = 0, ~~~~ T_1 = T^{0,0}, &~~ \Kc_t = \{ (1) \}, ~~ \forall t \in [0,T_1],
        \e*
        and we have the induction relationship
        \b*
                N_{T_{i+1}} ~=~ N_{T_i} + I_{i+1} - 1,&& T_{i+1} ~=~ T_i + \min_{k \in \Kc_{T_i}} T^{i, c(k)} ~=~ T_i + T^{i, c(K_{i+1})},
        \e*
        where $K_{i+1}$ denote the index of the particle which branches at time $T_{i+1}$. Let
        \b*
                \Kc_{T_{i+1}} ~:=~ \big \{ (K_{i+1},m) ~: 1 \le m \le I_{i+1} \big \} ~\cup~ \big \{ (k,1) ~: k \in \Kc_{T_i} \setminus \{ K_{i+1} \} \big \}.
        \e*
        In particular, if $I_{i+1} = 0$, then $\Kc_{T_{i+1}} = \big \{ (k,1) ~: k \in \Kc_{T_i} \setminus \{ K_{i+1} \} \big \}$.
        Clearly, at a branching time $T_i$, the particle $K_i$ branches into $I_i$ particles which are indexed by $(K_i,1), \cdots, (K_i,I_i)$, and all the other particles with index $k$ are re-indexed by $(k, 1)$. Let
        \b*
                N_t := N_{T_i} ~~\mbox{and}~~ \Kc_t := \Kc_{T_i}, &\mbox{for all}& t \in [T_i, T_{i+1}).
        \e*
        Then $(N_t)_{t \ge 0}$ is a continuous-time Markov process taking value in $\N$. Since it is possible that a particle dies with $k=0$ descendants, the branching system is subject to extinction in finite time horizon, i.e. $\P[N_t = 0~\mbox{for some}~t>0] > 0$. Furthermore, $(\Kc_t)_{t \ge 0}$ is a random process taking value in $ \cup_{n \in \N} 2^{\{1, \cdots, {n_0} \}^n}$, and $N_t = 0$ whenever $\Kc_t$ is empty.

        \begin{Example}
                We give an example of the branching birth-death process, with graphic illustration below, where $n_0 = 2$. The process starts with one particle indexed by $(1)$, and branches at time $T_1, \cdots, T_5$. The index of the branched particles are respectively $(1)$, $(1,1)$, $(1,2,1)$, $(1,1,2,1)$ and $(1,1,1,1,1)$.
                At terminal time $T$, the number of particles alive are $N_T = 5$ and
                \b*
                        \Kc_T &=& \big\{(1,1,1,1,1,1), (1,1,1,1,1,2), (1,1,2,1,1,1), (1,2,1,1,1,1), (1,2,1,2,1,1) \big\}.
                \e*
                \begin{itemize}
                    \item At time $T_1$, particle $(1)$ branches into two particles $(1,1)$ and $(1,2)$.
                    \item At time $T_2$, particle $(1,1)$ branches into $(1,1,1)$ and $(1,1,2)$, particle $(1,2)$ is reindexed by $(1,2,1)$.
                    \item At time $T_3$, particle $(1,2,1)$ branches into $(1,2,1,1)$ and $(1,2,1,2)$, the other two particles are reindexed by $(1,1,1,1)$ and $(1,1,2,1)$.
                    \item At time $T_4$, particle $(1,1,2,1)$ branches into one particle $(1,1,2,1,1)$, the other particles are reindexed.
                    \item At time $T_5$, particle $(1,1,1,1,1)$ branches into $(1,1,1,1,1,1)$ and $(1,1,1,1,1,2)$, the other particles are reindexed by $(1,1,2,1,1,1)$, $(1,2,1,1,1,1)$, $(1,2,1,2,1,1)$.
                \end{itemize}
        \end{Example}

                \setlength{\unitlength}{0.8cm}
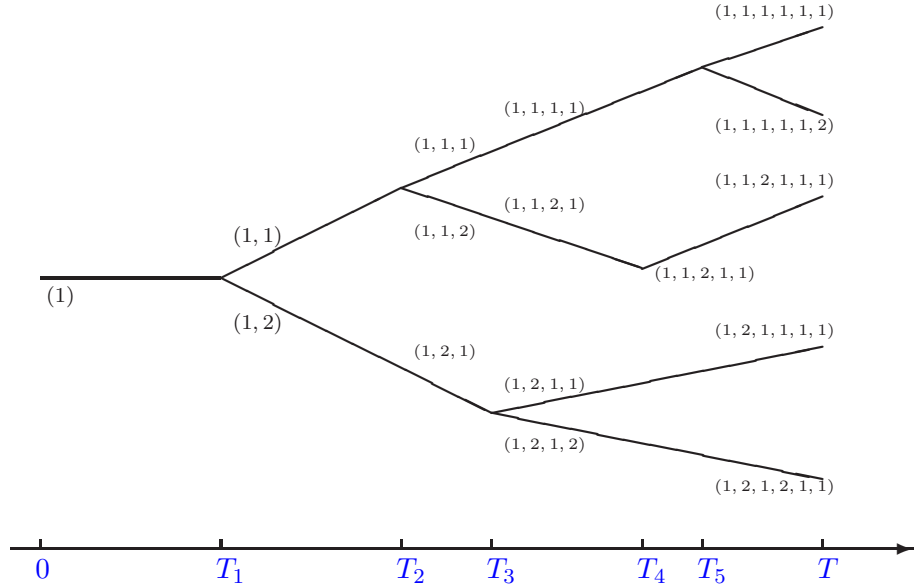
                \begin{picture}(20,10)
                        \thicklines

                        \put(0.5,0.5){\vector(1,0){15}}

                        \put(0.92,0){\cblue $0$}
                        \put(1,0.5){\line(0,1){0.15}}

                        \put(3.92,0){\cblue $T_1$}
                        \put(4,0.5){\line(0,1){0.15}}

                        \put(6.92,0){\cblue $T_2$}
                        \put(7,0.5){\line(0,1){0.15}}

                        \put(8.42,0){\cblue $T_3$}
                        \put(8.5, 0.5){\line(0,1){0.15}}

                        \put(10.92,0){\cblue $T_4$}
                        \put(11, 0.5){\line(0,1){0.15}}

                        \put(11.92,0){\cblue $T_5$}
                        \put(12, 0.5){\line(0,1){0.15}}

                        \put(13.92,0){\cblue $T$}
                        \put(14, 0.5){\line(0,1){0.15}}

                        \put(1, 5){\line(1,0){3}}

                        \put(4, 5){\line(2,1){3}}
                        \put(4, 5){\line(2,-1){4.5}}

                        \put(7, 6.5){\line(5,2){5}}
                        \put(7, 6.5){\line(3,-1){4}}

                        \put(8.5, 2.75){\line(5,1){5.5}}
                        \put(8.5, 2.75){\line(5,-1){5.5}}

                        \put(11, 5.15){\line(5,2){3}}

                        \put(12, 8.5){\line(3,1){2}}
                        \put(12, 8.5){\line(5,-2){2}}

                        \put(1.1, 4.6){{\scriptsize $(1)$}}

                        \put(4.2, 4.15){{\scriptsize $(1,2)$}}
                        \put(4.2, 5.6){{\scriptsize $(1,1)$}}

                        \put(7.2, 7.15){{\tiny $(1,1,1)$}}
                        \put(7.2, 5.7){{\tiny $(1,1,2)$}}
                        \put(7.2, 3.7){{\tiny $(1,2,1)$}}

                        \put(8.7, 7.75){{\tiny $(1,1,1,1)$}}
                        \put(8.7, 6.15){{\tiny $(1,1,2,1)$}}
                        \put(8.7, 3.15){{\tiny $(1,2,1,1)$}}
                        \put(8.7, 2.15){{\tiny $(1,2,1,2)$}}

                        \put(11.2, 5.){{\tiny $(1,1,2,1,1)$}}

                        \put(12.2, 9.35){{\tiny $(1,1,1,1,1,1)$}}
                        \put(12.2, 7.45){{\tiny $(1,1,1,1,1,2)$}}
                        \put(12.2, 6.55){{\tiny $(1,1,2,1,1,1)$}}
                        \put(12.2, 4.05){{\tiny $(1,2,1,1,1,1)$}}
                        \put(12.2, 1.45){{\tiny $(1,2,1,2,1,1)$}}

                \end{picture}

        \begin{Lemma}
                For every probability density sequence $(p_k)_{0 \le k \le {n_0}}$, we have $\lim_{n \to \infty} T_n = \infty$, a.s. In particular, the system is well defined from $0$ to $\infty$.
        \end{Lemma}
        \proof Without loss of generality, we can consider the case when $p_k = 0,~\forall k < n_0$ and $p_{n_0} = 1$.
        We first claim that $ N_t < \infty$ for all $t \ge 0$, it follows that $\sup\{n ~: T_n \le t \} < \infty$ for all $t \ge 0$ and hence $\lim_{n \to \infty} T_n = \infty$.
        Then to conclude, it is enough to prove that $N_t < \infty, ~ \forall t \ge 0$, which means that the population of the particles never explodes.
        It is then enough to use Example 2 of Kersting and Klebaner \cite{KerstingKlebaner} to finish the proof.
                \qed

        \vspace{2mm}

\paragraph{The branching Brownian motion}

        Suppose that in the same probability space $(\Om, \Fc, \P)$, there is a sequence of independent
        $d-$dimensional standard Brownian motions $(W^1, W^2, \cdots)$, which is also independent of the
        exponential random variables $(T^{i,j})_{i, j \ge 0}$ and multi-nomial random variables $(I_n)_{n \ge 1}$.
        We can then construct a branching Brownian motion which starts at time $t \ge 0$.

        For the first particle in the system indexed by $ k = (1)$, we associate it with a Brownian motion
         on $[t, \infty) $, defined by $B^{t,(1)}_{t+s} = W^1_s$, $\forall \; 0 \le s \le T_1$.
        Let $k = (k_1, \cdots, k_n) \in \Kc_{T_n}$ be the index of a living particle at time $T_n$, whose
        parent particle is indexed by $(k_1, \cdots, k_{n-1})$, we associate it with a Brownian motion
        between $[t, t+T_{n+1}]$, defined by
        \b*
            B^{t,k}_{t + s} &:=&
            \begin{cases}
                B^{t,(k_1, \cdots, k_{n-1})}_{t+s}, & \forall s \in [0, T_n], \\
                B^{t,(k_1,\cdots, k_{n-1})}_{t+T_n} ~+~ W_{s-  T_n}^{c(k)}, &\forall s \in [T_n, T_{n+1}].
            \end{cases}
        \e*
        By the strong Markov property of the Brownian motion, it is clear that conditioned on $(T^{i,j})_{i,j \ge 0}$ and $(I_n)_{n \ge 0}$, every process $(B^{t,k}_r )_{r \ge t}$ for $k \in \Kc_T$ is a Brownian motion.
        In particular, given two particles $k^1 = (k^1_1, \cdots, k^1_n)$ and $k^2 = (k^2_1, \cdots, k^2_n)$ such that $k^1_j = k^2_j$ for all $j=1, \cdots, i$, the associated Brownian motions $B^{t,k^1}$ and $B^{t,k^2}$ share the same path before time $t+T_i$.

\paragraph{The branching diffusion process}

        To construct a branching diffusion process, we first remark that for every $(t,\xb) \in \Lambda^0$, the SDE \eqref{eq:SDE} with initial condition $^{t,\xb}X_s = \xb_s,~ 0 \le s \le t$ has a
        unique strong solution $^{t,\xb}X$ adapted to the natural Brownian filtration, hence there
        is a progressively measurable function $\Phi^{t,\xb} : [t, T] \x C([t,T],\R^d) \to \R$ such
        that $^{t,\xb}X_s = \Phi^{t,\xb}(s,B_{\cdot}),~ \P_0-\mathrm{a.s.}$.

        Then a branching diffusion process $^{t,\xb}X^k$ is given by
        \be \label{eq:branch_diffu}
            ^{t,\xb}X_{t+s}^k &:=& \Phi^{t,\xb}(t+s, B^{t,k}_{\cdot}), ~~~\forall
            s \in \R^+ ~~\mbox{and}~~ k \in \Kc_s.
        \ee
				Moreover, for later uses, we extend $^{t,\xb}X^{(1)}$ on the whole interval $[0,T]$ by
				\b*
						^{t,\xb}X^{(1)}_s := \xb_s ~~ \forall s \le t &\mbox{and}& ^{t,\xb}X_s^{(1)} := \Phi^{t,\xb}(s, B^{t,(1)}_{\cdot}), ~~ \forall s \ge t.
				\e*

        \begin{Remark}
            By the flow property of the SDE \eqref{eq:SDE}, we have that for every $(t,\xb) \in \Lambda^0$, $r \le s$ and $k \in \Kc_s$,
            \be \label{eq:tower_prop_SDE}
                    \Phi^{t,\xb} \big( t+s, (B^{t,k}_u)_{u \ge t} \big) &=&
                    \Phi^{t+r, ~^{t, \xb}X^k} \big( t+s, (B^{t,k}_u)_{u \ge t+r} \big),~ \P-a.s.
            \ee
        \end{Remark}

        To conclude this subsection, we equip the above system with two filtrations.
        First, $\Fbb = (\Fcb_t)_{t \ge 0}$ with
        \b*
            \Fcb_t &:=& \sigma \big( (T_n, I_n, K_n) 1_{T_n \le t} + \partial 1_{T_n > t}, ~ n\ge 1 \big),
        \e*
        where $\partial$ denotes a cemetery point.
        Intuitively, $\Fbb$ is the filtration generated by the birth-death process. In
        particular, $T_n$ is a $\Fbb-$stopping time and $\Fcb_{T_n} = \sigma \big( (T_k, I_k, K_k)_{1 \le k \le n} \big) $.
        Next, for every $t \ge 0$, let $\Fbb^t = (\Fcb^t_{t+s})_{s \ge 0}$ be the filtration on the probability space $(\Om, \Fc, \P)$ generated by the branching diffusion process, which is defined by
        \be \label{eq:Fbb}
            \Fcb^t_{t+s} &:=& \Fcb_s ~\bigvee~ \sigma \big( \Kc_r,~ B^{(1)}_r, ~B^{t,k}_{t+r},~ 0 \le r \le s,~ k \in \Kc_s \big).
        \ee

\subsection{Branching diffusion representation of backward SDE}

        Using the branching diffusion process defined above, we can provide a representation of the solution to the decoupled FBSDE \eqref{eq:FBSDE_poly}.

        Let $(t,\xb) \in \Lambda^0$, we consider the branching diffusion process $(^{t,\xb}X^k)_{k \in \Kc_T}$ on $[t,T]$ defined in \eqref{eq:branch_diffu},
        where the probability sequence $p = (p_k)_{0 \le k \le n_0}$ satisfies that $p_k > 0$ whenever $|a_k|_0 \neq 0$.
        Denote
        \b*
            \Wc_{t,\xb}
            :=
            \Pi_{n =1}^{M_{T-t}} \Big( \frac{a_{I_n}(t+T_n, ~^{t,\xb}X^{K_n}_{\cdot})}{p_{I_n}} \Big),
            &\mbox{where}&
            M_{T-t}
            :=
            \sup\{n: t+T_n\le T\},
        \e*
        is the number of branchings occurred in the particles system between $t$ and $T$, with the convention that $\Pi_{n = 1}^0 := 1$.
        Our main representation formula is the following function:
        \be \label{eq:def_v}
            v(t,\xb)
            ~:=~
            \E^{\P} ~ \big[ \Psi_{t,\xb} \big]
            ~~&\mbox{with}&
            \Psi_{t,\xb}
            ~:=~
            \Wc_{t,\xb} ~\Pi_{k \in \Kc_{T-t}} ~ \psi\big(^{t, \xb}X_{\cdot}^k \big),
        \ee
where the integrability of $\Psi_{t,\xb}$ is verified in the
following result.

        \begin{Proposition} \label{prop:property_v}
            Suppose that Assumption \ref{assum:function_p} holds true. Then for every $(t,\xb) \in \Lambda^0$, the random variable $\Psi_{t,\xb}$ given in \eqref{eq:def_v} is integrable and the value function $v$ is uniformly bounded.
            Moreover, for every $M > 0$, there is a constant $C$ such that
            \b*
                \big|v(t, \om) - v(t', \om') \big|
                &\le&
                C \big(\sqrt{|t - t'|} + \| \om_{t \land \cdot} - \om'_{t'\land \cdot} \|_T  \big),
            \e*
            whenever $|(t,\om)| \le M$ and $|(t', \om')| \le M$.
        \end{Proposition}

The proof of Proposition \ref{prop:property_v} will be completed later in Section \ref{sec:proof_property_v}.

Our main result of the paper is the following representation
theorem. Let $^{0,\0}X$ be the unique strong solution to the SDE
\eqref{eq:SDE} in the probability space $(\Om, \Fc, \P)$, denote
        \be \label{eqd:Yb}
            Y^0_t &:=& v(t,^{0,\0}X_{\cdot}).
        \ee
        We also consider the BSDE \eqref{eq:FBSDE_poly} with generator
        \be \label{eq:F_poly}
                F_{n_0}(t, \xb, y) &:=& \beta \Big( \sum_{k=0}^{n_0} a_k(t,\xb) y^k ~-~ y \Big).
        \ee
        We define $\ell_{n_0}$ by
        \b*
                \ell_{n_0}(s)
                &:=&
                \beta \Big( \sum_{k=0}^{n_0} |a_k|_0 |\psi|_0^{k-1} s^k - s \Big), ~~ \forall s \ge 0.
        \e*
        It is clear that when Assumption \ref{assum:function_p} holds true
        for $\ell$, then $\ell_{n_0}$ satisfies also Assumption \ref{assum:function_p}. It
        follows from Proposition \ref{prop:Lip_BSDE} that the BSDE \eqref{eq:FBSDE_poly} with generator
        $F_{n_0}$ has a unique solution, denoted by $(Y,Z)$, such that $Y$ is uniformly bounded.

        \begin{Theorem} \label{theo:Branching_BSDE}
            Suppose that Assumption \ref{assum:function_p} holds true, and $(Y,Z)$ is the unique solution of BSDE \eqref{eq:FBSDE_poly} with generator $F_{n_0}$ (defined by \eqref{eq:F_poly}) such that $Y$ is uniformly bounded by $R_0$, the constant introduced in Lemma \ref{lemm:solution_ODE}.
            Then $Y^0 = Y$, $\P_0-$a.s.
        \end{Theorem}

The proof of this result will be provided in Section
\ref{sec:proof_mainthm} using the notion of viscosity solutions to a
path dependent PDE.

    \begin{Remark}
        The results in Proposition \ref{prop:property_v} and Theorem \ref{theo:Branching_BSDE} hold true for any probability sequence $p = (p_k)_{0 \le k \le n_0}$ satisfying that $p_k > 0$ whenever $|a_k|_0 \neq 0$. This implies that the integrability and expectation of $\Psi_{t,\xb}$ is independent of the choice of $p$.
        However, the variance of $\Psi_{t,\xb}$ does depend on $p$, where an upper bound is given by
        \be \label{eq:var_upperbound}
            \E \left[ \Pi_{n =1}^{M_{T-t}} \Big( \frac{|a_{I_n}|_0^2}{p_{I_n}^2} \Big) ~\Pi_{k \in \Kc_{T-t}} ~ |\psi|_0^2 \right].
        \ee
        Comparing $\Psi_{t,\xb}$ in \eqref{eq:def_v} with the integral part in \eqref{eq:var_upperbound}, it can be considered as a manipulation of the coefficients from $(a_k)_{k \ge 0}$ to $ \big(\frac{|a_k|_0^2}{p_k} \big)_{k \ge 0}$.
        Denote by $R^v$ the convergence radius of the sum $\sum_{k \ge 0} \frac{|a_k|_0^2}{p_k} s^k$,
        then using Proposition \ref{prop:property_v} with the new coefficients,
        the upper bound \eqref{eq:var_upperbound} is finite if and only if $\ell^v(1)$ one of the conditions $(\ell 1 - \ell 3)$ in Assumption \ref{assum:function_p},
        where
        \b*
            \ell^v(s) := \beta \Big( \sum_{k=0}^{n_0} \frac{|a_k|_0^2}{p_k} |\psi|_0^{2k-1} s^k -s \Big).
        \e*
    \end{Remark}

\subsection{Numerical algorithm by branching process}
\label{subsec:num_algo}

    The representation result in Theorem \ref{theo:Branching_BSDE} induces immediately a numerical algorithm
    for BSDE \eqref{eq:FBSDE_poly} by simulating the branching diffusion process. For numerical implementation, the branching times can be exactly simulated since they follow the exponential law, and
    the diffusion process can be simulated by a Euler scheme.

    Let $\Delta = (t_0, \cdots, t_n)$ be a discretization of the
    interval $[0,T]$, i.e. $0 = t_0 < \cdots < t_n = T$. Denote $|\Delta | := \max_{1 \le k \le n}
    (t_k - t_{k-1})$. To give the Euler scheme, we introduce the frozen coefficients $\mu^{\Delta}$ and $\sigma^{\Delta}$ by
    \b*
            \mu^{\Delta}(t, \xb) ~:=~ \mu(t_k, \xbh^{\Delta}) &\mbox{and}& \sigma^{\Delta}(t, \xb) ~:=~ \sigma(t_k, \xbh^{\Delta}), ~~~~ \forall t \in [t_k, t_{k+1}),
    \e*
    where $\xbh^{\Delta}$ denotes the linear interpolation of $(\xb_{t_0}, \cdots, \xb_{t_n})$ on the interval $[0,T]$. Then clearly the process $X^{\Delta}$ given by the SDE
    \be \label{eq:SDE_Delta}
            X^{\Delta}_t &=& \int_0^t \mu^{\Delta}(s, X^{\Delta}_{\cdot}) ds ~+~  \int_0^t \sigma^{\Delta}(s, X^{\Delta}_{\cdot}) dB_s,~~ \P_0-a.s.
    \ee
    can be simulated, which is also the Euler scheme of the SDE \eqref{eq:SDE}. By standard arguments
    using Gronwall's Lemma (see e.g. Kloeden and Platen \cite{KloedenPlaten} or Graham and Talay \cite{GrahamTalay} in the Markov case), we have the following error analysis result: Let $X$ be the solution process of \eqref{eq:SDE} with initial condition $(t,\xb) = (0,\0)$, $X^{\Delta}$ be the solution of \eqref{eq:SDE_Delta} and $\widehat{X}^{\Delta}$ denotes the linear interpolation of $(X^{\Delta}_{t_0}, \cdots, X^{\Delta}_{t_n})$ on $[0,T]$.

    \begin{Lemma} \label{lemm:error_Euler}
        There is a constant $C$ independent of the discretization $\Delta$ such that
        \b*
                    \E \Big[ \sup_{0 \le t \le T}
                    \Big( \big|X_t - X^{\Delta}_t \big|^2 + \big| X^{\Delta}_t - \widehat{X}^{\Delta}_t \big|^2 \Big)
                    \Big]
                    &\le&
                    C~ |\Delta|.
        \e*
    \end{Lemma}

    Moreover, for the BSDE \eqref{eq:FBSDE_poly} with a general generator function $F:[0,T] \x \Om^0 \x \R \to \R$ which admits a representation \eqref{eq:def_f}, we can approximate it by some polynomial $F_{n_0}$ of the form \eqref{eq:F_poly}. Let $F_{n_0}^{\Delta}(t,\xb,y) := \beta \big( \sum_{k=0}^{n_0} a^{\Delta}_k(t,\xb) y^k - y \big)$, where
    \b*
        a_k^{\Delta}(t,\xb) ~:=~ a_k(t_i, \xbh^{\Delta}) &\mbox{for every}& k = 0,\cdots, {n_0}~~\mbox{and}~~ t \in [t_i, t_{i+1}).
    \e*
    Further, under Assumption \ref{assum:function_p}, by simulating the branching diffusion process $(X^{\Delta, k})_{k \in \Kc_T}$, the numerical solution
    \b*
            Y^{\Delta}_0 &:=& \E \Big[ \Pi_{n=1}^{M_T} \Big( \frac{a^{\Delta}_{I_n}(T_n, \widehat{X}^{\Delta, K_n}_{\cdot})}{p_{I_n}} \Big) \Pi_{k \in \Kc_T} \psi \big(\widehat{X}^{\Delta, k}_{\cdot} \big) \Big]
    \e*
    is the solution of the following BSDE
    \b*
            Y_0 &=& \psi \big(\widehat{X}^{\Delta}_{\cdot} \big) ~+~
             \int_0^T F^{\Delta}_{n_0} \big( t, \widehat{X}^{\Delta}_{\cdot}, Y_t \big) dt ~- \int_0^T Z_t~ dB_t, ~~~ \P_0-a.s.
    \e*
    Finally, we provide an error estimation of the numerical solution:

    \begin{Proposition}
        Under Assumption \ref{assum:function_p}, there is a constant $C$ independent of $n_0$ and $\Delta$ such that
        \b*
            |Y^{\Delta}_0 - Y_0| &\le& C \big( |F - F^{\Delta}_{n_0}|_{L^{\infty}(\Lambda^0 \x [-R_0, R_0])} + \sqrt{\Delta} \big).
        \e*
    \end{Proposition}
    \proof This estimate follows from a direct application of the stability result of backward SDEs together with the error estimation in Lemma \ref{lemm:error_Euler},
    see Proposition 2.1 and their subsequent remark in El Karoui, Peng and Quenez \cite{ElKarouiBSDE}. \qed

    \begin{Remark}
        Let us consider an arbitrary Lipschitz generator $\overline F : \Lambda^0 \x \R$, such that the associated BSDE of the form \eqref{eq:FBSDE_poly} has a unique solution $(\overline Y, \overline Z)$ where $| \overline Y |$ is uniformly bounded by some constant $R_0 > 0$.
        One can also approximate the function $\overline F(y)$ by a polynomial function $\overline F^{n}(y) := \sum_{k=0}^n a^n_k y^k $.
        We may then conduct our analysis by formulating Assumption \ref{assum:function_p} on the coefficients $(a^n_k)_{0 \le k \le n}$ for all $n \ge 1$.
        The problem with this approach is that the convergence condition would depend on the approximating sequence of polynomials.
    \end{Remark}

\section{H\"{o}lder and Lipschitz regularity of $v$}
\label{sec:proof_property_v}

This section is devoted to the proof of Proposition
\ref{prop:property_v}. We first derive some estimates of the
birth-death process defined in Section
\ref{subsec:branch_diffusion}, then together with the tower
property, we can complete the proof of Proposition
\ref{prop:property_v}.

\subsection{Some estimates of the birth-death process}

We recall that $\Fbb = (\Fcb_t)_{0 \le t \le T}$ is the filtration
generated by the birth-death process defined in the end of Section
\ref{subsec:branch_diffusion}, and that the number of branchings
occurred in the system before time $t$ is   denote by $M_t ~:=~ \sup
\big \{ n ~: T_n \le t \big \}$. We also introduce:
    \b*
            \eta(t)
        &:=&
        \E^{\P} \Big[  \Big( \Pi_{n=1}^{M_t} \frac{|a_{I_n}|_0}{p_{I_n}} \Big) |\psi|_0^{N_t} \Big].
    \e*

    \begin{Lemma} \label{lemm:memoryless1}
        For every $0 \le s \le t$,
        \b*
            \E^{\P} \Big[ \Big( \Pi_{n=1}^{M_t} \frac{|a_{I_n}|_0}{p_{I_n}} \Big) |\psi|_0^{N_t} ~\Big |~ \Fcb_s
            \Big]
            &=&
            \Big( \Pi_{n=1}^{M_s} \frac{|a_{I_n}|_0}{p_{I_n}} \Big) (\eta(t-s))^{N_s},
        \e*
        and
        \b*
            \E^{\P} \Big[  \Big( \Pi_{n=1}^{M_t} \frac{|a_{I_n}|_0}{p_{I_n}} \Big) |\psi|_0^{N_t} 1_{T_1 \le t}
             ~\Big |~ \Fcb_{T_1} \Big]
            &=&
            \frac{|a_{I_1}|_0}{p_{I_1}} (\eta(t- T_1))^{I_1} 1_{T_1 \le t}.
        \e*
    \end{Lemma}
    \proof \rmi Let $Z$ be a random variable and $A \in \Fc$, then $\Lc^{\P}(Z)$ denotes the law of $Z$ and $\Lc^{\P}(Z|A)$ denotes the distribution of $Z$ conditioned on $A$ under the probability $\P$. We notice that for every $i, j \ge 1$ and $s > 0$,
        \b*
            \Lc^{\P}( T^{i,j} -s | T^{i,j} > s) &=& \Lc^{\P}(T^{i,j}) ~~=~~ \Ec(\beta).
        \e*
        Let $0 \le s \le t$, the law of number of branches between $s$ and $t$ (which equals to $M_t - M_s$) is completely determined by $N_s$, $(T^{i,j})_{i \ge M_s, j \ge 0}$ and $(I_i)_{i \ge M_s +1}$. It follows that
        \b*
                &&  \Lc^{\P} \big( M_t - M_s, ~(I_{M_s + i})_{ i \ge 1} ~\big|~ N_s = 1, ~k \in \Kc_s, ~ M_s = j \big) \\
                &=& \Lc^{\P} \big( M_t - M_s, ~(I_{M_s + i})_{ i \ge 1} ~\big|~ N_s = 1, ~k \in \Kc_s, ~ M_s = j, ~T^{M_s, c(k)} > s \big) \\
                &=& \Lc^{\P} \big( M_{t-s}, ~(I_i)_{i \ge 1} \big),
        \e*
        and hence
        \b*
            \Lc^{\P} \big( M_t - M_s, (I_{M_s + i})_{ i \ge 1} ~\big|~ N_s = 1 \big)
            &=&
            \Lc^{\P} \big( M_{t-s}, (I_i)_{i \ge 1} \big).
        \e*
        Since $N_t = N_s + \sum_{n = M_s +1}^{M_t} (I_n -1)$, we deduce that
        \b*
            \E^{\P} \Big[ \Big( \Pi_{n=M_s+1}^{M_t} \frac{|a_{I_n}|_0}{p_{I_n}} \Big) |\psi|_0^{N_t} ~\Big |~ N_s = 1 \Big]
            &=&
            \eta(t-s).
        \e*
        Moreover, since every particle branches independently to each other, we deduce that
        \b*
            \E^{\P} \Big[ \Big( \Pi_{n=M_s+1}^{M_t} \frac{|a_{I_n}|_0}{p_{I_n}} \Big) |\psi|_0^{N_t} ~\Big |~ N_s = i \Big]
            &=&
            (\eta(t-s))^i,
        \e*
        which implies that
        \b*
                \E^{\P} \Big[ \Big( \Pi_{n=M_s+1}^{M_t} \frac{|a_{I_n}|_0}{p_{I_n}} \Big) |\psi|_0^{N_t} ~\Big |~ N_s \Big]
                &=&
                (\eta(t-s))^{N_s}.
        \e*
        And hence
        \b*
                \E^{\P} \Big[ \Big( \Pi_{n=1}^{M_t} \frac{|a_{I_n}|_0}{p_{I_n}} \Big) |\psi|_0^{N_t} ~\Big |~ \Fcb_s \Big]
                &=&
                \E^{\P} \Big[ \Big( \Pi_{n=1}^{M_s} \frac{|a_{I_n}|_0}{p_{I_n}} \Big)
                \Big( \Pi_{n=M_s + 1}^{M_t} \frac{|a_{I_n}|_0}{p_{I_n}} \Big) |\psi|_0^{N_t} ~\Big |~ \Fcb_s \Big] \\
                &=&
                \Big( \Pi_{n=1}^{M_s} \frac{|a_{I_n}|_0}{p_{I_n}} \Big) (\eta(t-s))^{N_s}.
        \e*

        \rmii We next prove the second equality, we notice that $(I_i)_{i \ge 2}$ and $(T^{i,j})_{i \ge 2,~ j \ge 0}$ are all independent of $(T_1, I_1)$ under the probability $\P$. Let us consider a family of conditional probabilities $(\P_{s,i})_{s \in \R^+, i \in \{0, \cdots, n_0 \}}$ of $\P$ w.r.t. the $\sigma-$field generated by $(T_1, I_1)$. 
        Under every conditional probability $\P_{s,i}$, the law of $(M_t, N_t) 1_{s \le t}$ depends only on $(I_j)_{j \ge 2}$ and $(T^{j,l})_{j \ge 2,~ l \ge 0}$. Considering in particular $i =1$, we have
        \b*
                \Lc^{\P_{s,1}} \big( M_t - M_s,~ (I_i)_{i \ge 2} \big)
                &=&
                \Lc^{\P} \big( M_{t-s}, ~(I_i)_{i \ge 1} \big).
        \e*
        And hence
        \b*
            \E^{\P_{s,1}} \Big[  \Big( \Pi_{n=M_s + 1}^{M_t} \frac{|a_{I_n}|_0}{p_{I_n}} \Big) |\psi|_0^{N_t} 1_{s \le t} \Big]
            &=&
            \eta (t- s) 1_{s \le t}.
        \e*
        Moreover, by the independence of the evolution of $i$ particles under $\P_{s,i}$, we get
        \b*
            \E^{\P_{s,i}} \Big[  \Big( \Pi_{n=M_s + 1}^{M_t} \frac{|a_{I_n}|_0}{p_{I_n}} \Big) |\psi|_0^{N_t} 1_{s \le t} \Big]
            &=&
            (\eta (t- s))^i 1_{s \le t},
        \e*
        which implies that
        \b*
            \E^{\P} \Big[  \Big( \Pi_{n=1}^{M_t} \frac{|a_{I_n}|_0}{p_{I_n}} \Big) |\psi|_0^{N_t} 1_{T_1 \le t}
             ~\Big |~ \Fcb_{T_1} \Big]
            &=&
            \frac{|a_{I_1}|_0}{p_{I_1}}
            \E^{\P_{T_1,I_1}} \Big[  \Big( \Pi_{n=M_{T_1} + 1}^{M_t} \frac{|a_{I_n}|_0}{p_{I_n}} \Big) |\psi|_0^{N_t} 1_{T_1 \le t} \Big]\\
            &=&
            \frac{|a_{I_1}|_0}{p_{I_1}} (\eta(t- T_1))^{I_1} 1_{T_1 \le t},
        \e*
        since $M_{T_1} = 1$ by its definition.
        And we hence conclude the proof. \qed

  \begin{Lemma} \label{lemma:nonexplosion}
    Suppose that for some $t \ge 0$, $\eta(t) < \infty$. Then there is $ \delta > 0$ such that $\eta(s) < \infty$, for every $s \in [t, t+ \delta]$.
  \end{Lemma}
    \proof First, it follows from Lemma \ref{lemm:memoryless1} that for every $t, ~\delta \ge 0$,
            \b*
        \eta(t+\delta) &=& \E^{\P} \Big[  \Big( \Pi_{n=1}^{M_{\delta}} \frac{|a_{I_n}|_0}{p_{I_n}} \Big) (\eta(t))^{N_{\delta}} \Big].
      \e*
        Let us consider another pure birth process $(\Nt_t, \Kct_t)$ on a probability space $(\Omt, \Fct, \Pt)$ with the same constant characteristic $\beta$ and another probability sequence $(\pt_k)_{0 \le k \le {n_0}}$ such that $\pt_{n_0} = 1$. We suppose without loss of generality that ${n_0} \ge 2$ and denote $C := \max_{0 \le k \le {n_0}} \frac{|a_k|_0}{p_k} + \eta(t)$. Then clearly it is enough to prove that $\E ^{\Pt} \big[ C^{\Nt_{\delta}} \big] < \infty$ for some $\delta > 0$ to conclude the proof. The distribution of $\Nt_{\delta}$ can be computed explicitly (see e.g. Athreya and Ney \cite[Chapiter III.5, P109]{AthreyaNey}) and satisfies that for some constant $\tilde C > 0$,
        \b*
                \P\big[\Nt_{\delta} = n\big]
                ~\le~
                \tilde C \kappa_{\delta}^n
                &\mbox{with}&
                \kappa_{\delta}
                ~:=~
                \big(1 - e^{-\delta\beta(n_0-1)}\big)^{1/({n_0}-1)}.
        \e*
        Then for $\delta > 0$ small enough, $\kappa_{\delta}$ is small enough such that $\E ^{\Pt} \big[ C^{\Nt_{\delta}} \big] < \infty$.
        \qed

    \vspace{2mm}

    The birth-death system is closely related to ODE \eqref{eq:ODE}. Let us define
    \be \label{eq:Branching_Integrability}
        D_t
        &:=&
      \E^{\P} \Big[ (1 \vee M_t) (1 \vee N_t) \Big( \Pi_{n=1}^{M_t} \frac{|a_{I_n}|_0}{p_{I_n}} \Big) |\psi|_0^{N_t} \Big].
    \ee

    \begin{Proposition} \label{prop:Cvg_GW}
        Suppose that Assumption \ref{assum:function_p} holds true. Then
        \be\label{eq:Cvg_GW}
            \sup_{0 \le t \le T} \eta(t) < \infty
            &\mbox{and}&
            \sup_{0 \le t \le T}
            D_t
            < \infty.
        \ee
    \end{Proposition}
    \proof We first observe that $\eta(0) = |\psi|_0$ by its definition, and it follows from Lemma \ref{lemm:memoryless1} that
        \b*
            \eta(t) &=& \E[ |\psi|_0 1_{T_1 > t}] ~+~  \E\Big[  \sum_{k=0}^{n_0} |a_k|_0 (\eta(t-T_1))^k 1_{T_1 \le t} \Big] \\
            &=&
            \eta(0) e^{-\beta t} ~+~ \int_0^t \beta e^{-\beta s} \Big( \sum_{k=0}^{n_0} |a_k|_0 (\eta(t-s))^k \Big) ds \\
            &=&
            e^{-\beta t} \Big( \eta(0)  ~+~ \int_0^t \beta e^{\beta s} \Big( \sum_{k=0}^{n_0} |a_k|_0 (\eta(s))^k \Big) ds \Big).
        \e*
        Suppose that $T_0 := \inf \{ s : \eta(s) = \infty \} \le T$, then it follows from Lemma \ref{lemm:solution_ODE} and Remark \ref{rem:eqv_ODE} that $\eta(t) = \rho(t)|\psi|_0, ~\forall t \in [0, T_0)$, where $\rho$ is the unique solution to ODE \eqref{eq:ODE}.
        Therefore, it follows still by Lemma \ref{lemm:solution_ODE} and Remark \ref{rem:eqv_ODE} that $\eta(T_0) = \rho(T_0) |\psi|_0 < \infty$, and hence $\eta(t) < \infty, ~ \forall t \in [0, T_0 + \delta]$ for some constant $0 < \delta < T - T_0$ by Lemma \ref{lemma:nonexplosion}.
        This contradicts the definition of $T_0$, and hence $T_0 > T$ and $\eta(T) < \infty$.
        Since $\eta$ is increasing, this provides the first claim in \eqref{eq:Cvg_GW}.

    We next denote
    \b*
        \eta_{\eps}(t) &:=&  \E^{\P} \Big[  \Big( \Pi_{n = 1}^{M_t} \frac{|(1+\eps)a_{I_n}|_0}{p_{I_n}} \Big) |(1+\eps)\psi|_0^{N_t} \Big].
    \e*
        In spirit of Remark \ref{rem:ODE_stability}, we know that for $\eps > 0$ small enough, $\eta_{\eps}(T) < \infty$. It follows that
    \b*
        \sup_{0 \le t \le T} D_t
        ~:=~
        \sup_{0 \le t \le T}
        \E^{\P} \Big[  (1 \vee M_t) (1 \vee N_t) ~ \Big( \Pi_{n, ~ T_n \le t} \frac{|a_{I_n}|_0}{p_{I_n}} \Big) |\psi|_0^{N_t} \Big]
        &<& \infty,
    \e*
    since there is some constant $C_{\eps}>0$ such that $n < C_{\eps} (1+\eps)^n$ for every $n \ge 0$. And we hence conclude the proof.    \qed

\subsection{Proof of Proposition \ref{prop:property_v}}
\label{subsec:proof_property_v}

        In preparation of the proof, we first provide a tower property of the branching diffusion process. Let $(t,\xb) \in \Lambda^0$ and $\tau : \Om^0 \to \R^+$ be a $\F^0-$stopping time such that $\tau \ge t$, then $\tauh := \tau( ~^{t,\xb}X^{(1)}_{\cdot})$ is a $\Fbb^t-$stopping time in the probability space $(\Om, \Fc, \P)$, which is clearly independent of $T^1$.

    \begin{Lemma} \label{lemm:memoryless2}
        Suppose that Assumption \ref{assum:function_p} holds true, let $(t,\xb) \in \Lambda^0$, $0 \le s \le T-t$ and $\tauh$ be given above.
        Then we have
        \be \label{eq:tower_prop_branching_diffusion1}
            \E^{\P} \big[ \Psi_{t,\xb} ~\big |~ \Fcb^t_{t+s} \big]
            &=&
            \Big( \Pi_{n=1}^{M_s} \frac{a_{I_n}(t+T_n, ~^{t,\xb}X^{K_n})}{p_{I_n}} \Big) ~\Pi_{k \in \Kc_s} v(t+s, ~^{t,\xb}X^k_{\cdot})
        \ee
        and
        \be \label{eq:DPP}
            v(t,\xb)
            &=&
            \E^{\P} \Big[ v(\tauh, ~^{t,\xb}X_{\cdot}^{(1)}) 1_{t+T_1 > \tauh} \nonumber \\
            && ~~~~+~
            \frac{a_{I_1} \big( t+T_1, ~^{t,\xb}X^{(1)}_{\cdot} \big)}{p_{I_1}} v^{I_1} \big( t+T_1,~ ^{t,\xb}X^{(1)}_{\cdot} \big) 1_{t+T_1 \le \tauh} \Big].
        \ee
    \end{Lemma}
        \proof First, following the arguments of Lemma \ref{lemm:memoryless1}, we know
        \b*
                && \Lc^{\P} \big(M_t - M_s, ~(I_{M_s + j})_{j \ge 1}, ~(W^{M_s + l})_{l \ge 1} ~\big|~ N_s = 1, ~M_s = i \big) \\
                &=&
                \Lc^{\P} \big( M_{t-s}, ~(I_j)_{j \ge 1}, ~(W^l)_{l \ge 1} \big).
        \e*
        Together with the flow property of SDE in \eqref{eq:tower_prop_SDE}, it follows that
        \b*
            && \E^{\P} \Big[ \Big( \Pi_{n = M_s+1}^{M_{T-t}} \frac{a_{I_n}(t+T_n, ~^{t,\xb}X^{K_n}_{\cdot})}{p_{I_n}} \Big)
                \Pi_{k \in \Kc_{T-t}} \psi\big(~^{t, \xb} X_{\cdot}^{k} \big)\\
            && ~~~~~~~~~~~~~~~~~~~~ \Big|~ N_s = 1,~ M_s = i,~ k \in \Kc_s, ~ (^{t,\xb}X_r^k)_{t \le r \le t+s} \Big]
            ~=~
            v(t+s, ~^{t,\xb}X_{\cdot}^k).
        \e*
        Then by the independence of evolution of every particle in $\Kc_s$, \eqref{eq:tower_prop_branching_diffusion1} holds true.

        For the second equality, we consider a regular conditional probability distribution (r.c.p.d.) $(\P_{\omh})_{\omh \in \Om}$ of $\P$ w.r.t. $ \sigma( B^{(1)}_{\tauh \land \cdot}) $ (see also Stroock Varadhan \cite{SV} for the notion of r.c.p.d.).
        Then for every $\omh \in \Om$, we have $\P_{\omh} \Big(B^{(1)}_s = B^{(1)}_s(\omh),~ 0 \le s \le \tauh(\omh) \Big) = 1$ and
        $\Big(B^{(1)}_s, ~ s \ge \tauh(\omh) \Big)$ is still a standard Brownian motion under $\P_{\omh}$.
        In particular, $\P_{\omh} \Big(X^{(1)}_s = X^{(1)}_s(\omh),~ 0 \le s \le \tauh(\omh) \Big) = 1$.
        Further, since $T_1 = T^{0,0}$ is independent of the Brownian motions $B^{(1)}$, then $T^{0,0}$ is still an exponential random variable under $\P_{\omh}$, and
        \b*
            \Lc^{\P_{\omh}} \Big(T^{0,0} - (\tauh(\omh) - t) ~\Big|~ T^{0,0} > (\tauh(\omh) - t) \Big)
            &=&
            \Lc^{\P} \big( T^{0,0} \big)
            ~~=~~
            \Ec(\beta).
        \e*
        By adding $\tauh(\omh)$ on each side, it follows that
        \b*
            \Lc^{\P_{\omh}} \Big( t+ T^{0,0} ~\Big|~ t + T^{0,0} > \tauh(\omh) \Big)
            &=&
            \Lc^{\P} \big( \tauh(\omh) + T^{0,0} \big).
        \e*
        By the expression of $\Psi_{t,\xb}$ and the fact that $\P_{\omh} \Big(X^{(1)}_s = X^{(1)}_s(\omh),~ 0 \le s \le \tauh(\omh) \Big) = 1$, we then have
        \b*
            \Lc^{\P_{\omh}} \Big( \Psi_{t,\xb} ~\Big|~ t + T^{0,0} > \tauh(\omh) \Big)
            &=&
            \Lc^{\P} \Big( \Psi_{\tauh(\omh), X^{(1)}_{\cdot}(\omh)} \Big).
        \e*
        Taking expectations, it follows that
        \b*
            \E^{\P_{\omh}} \big[ \Psi_{t,\xb} ~1_{t+T_1 > \tauh(\omh)} ~\big|~ t+T_1 > \tauh(\omh) \big]
            &=&
            \E^{\P_{\omh}} \big[ \Psi_{t,\xb} ~\big|~ t+T_1 > \tauh(\omh) \big]\\
            &=&
            v \big(\tauh(\omh), ~^{t,\xb} X_{\cdot}^{(1)}(\omh) \big),
        \e*
        and hence by the independence of $T_1$ to $~^{t,\xb}X^{(1)}$ and $\tauh$, we have
        \b*
            \E^{\P} \big[ \Psi_{t,\xb} ~1_{t+ T_1 > \tauh} \big]
            &=&
            \E^{\P} \big[ \Psi_{t,\xb} ~1_{t+ T_1 > \tauh} \big| t + T_1 > \tauh \big]
            ~ \P(t + T_1 > \tauh)\\
            &=&
            \E^{\P} \Big[ v \big(\tauh, ~^{t,\xb} X_{\cdot}^{(1)} \big) \Big]
            ~ \P(t + T_1 > \tauh)
            ~=~
            \E^{\P} \Big[ v \big(\tauh, ~^{t,\xb} X_{\cdot}^{(1)} \big) 1_{t+ T_1 > \tauh} \Big].
        \e*
        Further, using similar arguments as in Lemma \ref{lemm:memoryless1}, by considering the distribution of $\Psi_{t,\xb} 1 _{t + T_1 \le \tauh}$ conditioned on $\Fcb^1_{T_1}$, we get
        \b*
            \E^{\P} \big[ \Psi_{t,\xb} ~1_{t+ T_1 \le \tauh} \big]
            &=&
            \E^{\P} \Big[
            \frac{a_{I_1} \big( t+T_1, ~^{t,\xb}X^{(1)}_{\cdot} \big)}{p_{I_1}} v^{I_1} \big( t+T_1, ~^{t,\xb}X^{(1)}_{\cdot} \big) 1_{t+T_1 \le \tauh} \Big],
        \e*
        which concludes the proof. \qed

        \vspace{2mm}

        \noindent {\bf Proof of Proposition \ref{prop:property_v}}.
				\rmi First, it follows immediately from Proposition \ref{prop:Cvg_GW} that $\Psi_{t,\xb}$ is integrable and $|v(t,\xb)| \le \rho(T-t) |\psi|_0 \le R_0$.\\
     \rmii Let $t \in [0,T]$ and $\xb_1, ~ \xb_2 \in \Om^0$.
     It follows then by Lemma \ref{lemm:SDE_Lip} together with Cauchy-Schwartz inequality, that for every $s \in [t,T]$ and $k \in \Kc_s$:
    \b*
        \E^{\P} \Big[ \sup_{t \le r \le s} \big| ^{t,\xb_1}X^{k}_r - ~^{t,\xb_2}X^k_r \big| \Big]
        &\le&
        C \big( 1 + \|\xb_1\|_t + \| \xb_2 \|_t \big) ~ \| \xb_1 - \xb_2 \|_t,
    \e*
    for some constant $C$ independent of $\xb_1,~ \xb_2$. Then using the fact that $(a_k)_{0 \le k \le {n_0}}$ and $\psi$ are all Lipschitz in $\xb$,
        \b*
                \big| v(t,\xb_1) - v(t,\xb_2) \big|
                &\le&
                \E^{\P} \big[ \big| \Psi_{t,\xb_1} - \Psi_{t,\xb_2} \big| \big]
                ~\le~
                C ~D_t~ \E^{\P} \big[ \big\| ~^{t,\xb_1}X - ~^{t,\xb_2}X \big\|_T \big] \\
                &\le&
                C \big( 1 + \|\xb_1\|_t + \| \xb_2 \|_t \big) ~ \| \xb_1 - \xb_2 \|_t,
        \e*
        where $D_t$ is defined in \eqref{eq:Branching_Integrability}.\\
    \rmiii  Let $0 \le s \le t \le T$, then it follows from Lemma \ref{lemm:memoryless2} that
        \b*
                \big| v(s, \xb) - v(t,\xb_{s \land \cdot}) |
                &\le&
                \Big| \E^{\P} \Big[ \Pi_{n=1}^{M_{t-s}} \frac{a_{I_n}(T_n, ~^{s,\xb}X^{K_n}_{\cdot})}{p_{I_n}} ~\Pi_{k \in \Kc_{t-s}} v(t, ~^{s,\xb}X^k_{\cdot}) \Big] - v(t,\xb_{s \land \cdot}) \Big| \\
                &\le&
                C \big( \sup_{s \le r \le t} D_r \big) \E^{\P} \big[ \sup_{r \in [s,t]} \big| \xb_s - ^{s,\xb}X^k_r \big| \big] \\
                && +~
                \Big| \E^{\P} \Big[ \Pi_{n=1}^{M_{t-s}} \frac{a_{I_n}(t, \xb_{s \land \cdot})}{p_{I_n}} ~ \big( v(t, \xb_{s \land \cdot}) \big)^{N_{t-s}} \Big] - v(t,\xb_{s \land \cdot})  \Big| \\
                &\le&
                C (1 + \| \xb \|_s) \sqrt{t-s} ~+~ | \phi(t) - \phi(s) |,
        \e*
        where $\phi$ is the unique solution of the ODE
        \b*
                \phi'(r) = \beta \Big( \sum_{k=0}^{n_0} a_k(t,\xb) \phi^k(r) - \phi(r) \Big) &\mbox{with terminal condition}& \phi(t) = v(t,\xb_{s \land \cdot}).
        \e*
        Moreover, by comparison principle of ODE, $|\phi(r)| \le \rho(r), ~ \forall r \in [s, t]$. Then $| \phi(t) - \phi(s)| \le C(t-s)$ for some constant $C$ independent of $(s,t,\xb)$, which implies that $v$ is locally $(1/2)-$H\"older in $t$.
     \qed

        \begin{Remark}
            When $(a_k)_{0 \le k \le {n_0}}$ and $\psi$ are all constants, the value function $v(t,\xb)$ is independent of $\xb$ and $t \mapsto v(T-t,\xb) |\psi|^{-1}_0$ is a solution to ODE \eqref{eq:ODE}. Therefore, in spirit of Lemma \ref{lemm:solution_ODE}, Assumption \ref{assum:function_p} is also a necessary condition
            for the integrability of $\Psi_{0,0}$.
        \end{Remark}

\section{The branching diffusion representation result}
\label{sec:proof_mainthm}

This section is devoted to the proof of Theorem
\ref{theo:Branching_BSDE}.

We first consider a class of semi-linear parabolic path-dependent
PDEs (PPDEs) and introduce a notion of viscosity solution, following
Ekren, Keller, Touzi and Zhang \cite{EkrenKellerTouziZhang} and
Ekren, Touzi and Zhang \cite{EkrenTouziZhang1, EkrenTouziZhang2}.
Our objective is to show that the value function $v$, defined by our
branching diffusion representation, and the $Y-$component of the
BSDE are viscosity solutions of the same path-dependent PDE. Then,
our main result follows from a uniqueness argument.

\subsection{Viscosity solutions of PPDEs and FBSDEs}
\label{subsec:PPDE}

We consider a PPDE which is linear in the first and second order
derivatives of the solution function. This is a simpler context than
that of \cite{EkrenKellerTouziZhang, EkrenTouziZhang1,
EkrenTouziZhang2}. As a consequence, following Remark 3.9 in
\cite{EkrenTouziZhang1}, we use a simpler definition of viscosity
solutions. We shall also provide an (easy) adaptation of the
arguments in \cite{EkrenTouziZhang1} which relaxes their boundedness
conditions, thus allowing the terminal condition and the generator
to have linear growth.

\subsubsection{Differentiability on the canonical space}

For all $t \in [0,T]$, we denote by $\Om^t := \{ \om \in C([t,T],
\R^d) : \om_t = 0 \}$ the shifted canonical space, $B^t$ the shifted
canonical process on $\Om^t$, $\F^t$ the shifted canonical
filtration generated by $B^t$, $\P^t_0$ the Wiener measure on
$\Om^t$ and $\Lambda^t := [t,T] \x \Om^t$.

        For $s \le t$, $\om \in \Om^s$ and $\om' \in \Om^t$, define the concatenation path $\om \ox_t \om' \in \Om^s$ by
        \b*
                (\om \ox_t \om')(r) &:=& \om_r 1_{s \le r < t} ~+~ (\om_t + \om'_r)1_{t \le r \le T}, ~~\forall r \in [s,T].
        \e*
        Let $\xi \in \Fc^0_T$ and $V$ be a $\F^0-$progressive process, then for every $(t,\om) \in \Lambda^0$, we define $\xi^{t,\om} \in \Fc^t_T$ and $(V^{t,\om}_s)_{t \le s \le T}$ by
        \be \label{eq:shifted_proc}
                \xi^{t,\om}(\om') ~:=~ \xi(\om \ox_t \om'),
                &&
                V_s^{t,\om}(\om') ~:=~ V_s(\om \ox_t \om'), ~~ \forall \om' \in \Om^t.
        \ee

Following Ekren, Touzi and Zhang \cite{EkrenTouziZhang1, EkrenTouziZhang2}, we
define some classes of processes in $\Lambda^t$, $t\ge 0$. Let
$C^0(\Lambda^t)$ be the collection of all $\F^t-$progressive
processes which are continuous under the norm $d_{\infty}$, where
\b*
        d_{\infty} \big( (s,\om), (s', \om') \big)
        &:=&
        |s-s'| ~+~ \sup_{t \le r \le T} | \om_{s \land r} - \om'_{s' \land r} |,
        ~~~\forall (s,\om), (s',\om') \in \Lambda^t.
\e*
Denote by $C^0_b(\Lambda^t)$(resp. $UC(\Lambda^t)$) the
collection of functions in $C^0(\Lambda^t)$ which are uniformly
bounded (resp. uniformly continuous), and $UC_b(\Lambda^t) :=
UC(\Lambda^t) \cap C_b^0(\Lambda^t)$.

Next, denote by $X^{0, t,\xb}$ the solution of the SDE on $(\Om^t, \Fc^t_T, \P_0^t)$:
    \be \label{eq:SDEt}
    		X_s = \xb_s, ~\forall s \le t &\mbox{and}& X_s = \xb_t + \int_t^s \mu(r, X_{\cdot}) dr + \int_t^s \sigma(r, X_{\cdot}) dB^t_r,~\forall s > t.
    \ee
    Clearly, $X^{0,t,\xb}$ under $\P_0^t$ has the same law as that of $^{t,\xb}X$ introduced in \eqref{eq:SDE} under $\P_0$.
We denote the induced measure on the shifted space $\Omega^t$ by:
 \be \label{eq:Pt}
 \P_{t,\xb}
 :=
 \P^t_0 \circ \big(X^{0,t,\xb}-\xb_t\big)^{-1}
 &\mbox{and}&
 \P_X:=\P_{0,\mathbf{0}}.
 \ee

\begin{Remark} \label{rem:rcpd}
	Let $(t,\xb) \in \Lambda^0$, $\tau \ge t$ be a $\F^t-$stopping time on $\Om^t$, $\xi \in \Fc^t_T$ and
	$(\P_{\om})_{\om \in \Om}$ be a regular conditional probability distribution (r.c.p.d., see Stroock-Varadhan \cite{SV})
	of $\P_{t,\xb}$ w.r.t. $\Fc^0_{\tau}$, then clearly, $\E^{\P_{\om}}[ \xi] = \E^{\P_{\tau(\om), \om}}[ \xi^{\tau(\om), \om}]$ for $\P_{t,\xb}-$a.s. $\om \in \Om$.
\end{Remark}

For every $s \in [0,T)$ and $u:\Lambda^s\longrightarrow\R$, we
introduce the Dupire \cite{Dupire} right time-derivative of $u$
defined by the following limit, if exists,
    $$
        \partial_t u(t,\om)
        :=
        \lim_{h \downarrow 0} \frac{u(t+h,\om_{\cdot \land t})
                                    - u(t,\om)}
                                   {h},
        ~t<T,~\mbox{and}~
    \partial_t u(T,\om)
                :=
                \lim_{t <T, t \to T} \partial_t u(t,\om).
    $$

    \begin{Definition} \label{def:smooth_functions}
        Let $u$ be a process $C^0(\Lambda^t)$. We say $u \in C^{1,2}(\Lambda^t)$ if $\partial_t u \in C^0(\Lambda^t)$ and there exist $\partial_{\om} u \in C^0(\Lambda^t, \R^d)$, $\partial_{\om \om}^2 u \in C^0(\Lambda^t, \S^d)$ such that for all $r \ge t$,
        \be \label{eq:fct_Ito}
        d u_r
        &=&
        (\partial_t u)_r dr
        + (\partial_{\om} u)_r \cdot d B_r
        + \frac{1}{2} (\partial_{\om \om} u)_r
          : d \langle B \rangle_r,
        ~~\P_{t,\xb}-a.s.
        \ee
        If, in addition, $u \in C^0_b(\Lambda^t)$, we then say $u \in C^{1,2}_b(\Lambda^t)$.
    \end{Definition}

It is clear, for $s \le t$, $\om \in \Om^0$ and $u \in
C^{1,2}(\Lambda^s)$, we have $u^{t,\om} \in C^{1,2}(\Lambda^t)$.

    Finally, for all $ t \in [0,T]$, we denote by $\Tc^t$ the collection of all $\F^t-$stopping times $\tau$
    such that $\{\om ~: \tau(\om) > s \}$ is an open set in $(\Om^t, \| \cdot \|_T)$ for all $s \in [t,T]$, and by $\Tc^t_+$ the collection of stopping times $\tau \in \Tc^t$ such that $\tau > t$.
The set $\Lambda^t(\tau) := \{ (t,\om) \in \Lambda^t ~: t <
\tau(\om) \}$ is the corresponding localized canonical space, and we
define similarly the spaces $C^0(\Lambda^t(\tau) )$,
$C^{1,2}(\Lambda^t(\tau) )$, etc.

\subsubsection{A path-dependent PDE}

In this section, we do not need the restriction that the generator has a power series representation in $y$ as in \eqref{eq:def_f}. We
then consider a slightly more general generator $\widehat{F}:
\Lambda^0 \x \R \to \R$ such that $(t,\om) \longmapsto
\widehat{F}(t,\om,y)$ is $\F^0-$progressive for every $y \in \R$.
Consider the second order path-dependent differential operator:
    \be \label{eq:Lc}
        \Lc \varphi
        &:=&
        \partial_t \varphi
        ~+~ \mu \cdot \partial_{\om} \varphi
        ~+~ \frac{1}{2}\sigma \sigma^T:\partial_{\om \om}^2\varphi.
    \ee
Given a $\Fc_T-$measurable r.v. $\xi:\Omega^0\longrightarrow\R$, we
consider the path-dependent PDE:
    \be \label{eq:PPDE}
        -~\big\{ \Lc u + \widehat{F}(\cdot, u) \big\} (t,\om)
        ~=~
        0,
        &&\forall (t,\om) \in [0,T) \x \Om^0,
    \ee
    with terminal condition $u(T,\omega) = \xi(\omega), ~\forall \om \in \Om^0$.

    \begin{Assumption} \label{assum:fg}
    There is a constant $C$ such that $\sup_{t\le T} | \widehat{F}(t,\mathbf{0},0)|\le C$, and
            \b*
                \big|\widehat{F}(t, \om, y)
                     -\widehat{F}(t, \om', y') \big|
                +|\xi(\om) - \xi(\om')|
                &\le&
                C \big( |y - y'|  +  \| \om - \om' \|_T \big),
            \e*
for every $t \in [0,T]$, $(\om, y), ~(\om', y') \in \Om^0 \x \R$.
    \end{Assumption}

We denote by $\Ucb$ the class of functions $u$ defined on $\Lambda^0$ satisfying, for every $M > 0$, there is some continuity modulus $\rho_M$ such that
\b*
		u(t,\om) - u(t', \om') \le \rho_M \big(d_{\infty}\big( (t,\om), (t', \om') \big) \big),
		~\mbox{whenever}~t \le t' ~\mbox{and}~ \|\om \| \le M, \| \om'\| \le M,
\e*
and by $\Ucu$ the class of functions $u$ such that $-u \in \Ucb$;
we next introduce, for every $\F^0-$adapted process $u$, two classes of test functions:
     $$\begin{array}{cc}
        \Acb u(t,\om)
        :=
        \Big\{ \varphi \in C^{1,2}(\Lambda^t):
               \exists\;\H \in \Tc^t_+,~
               (\varphi- u^{t,\omega})_t(\0)
               =
               \displaystyle\min_{\tau \in \Tc^t}
               \E^{\P_{t,\om}}
               \big[ (\varphi - u^{t,\om})_{\tau \land \H} \big]
               \Big\},
               \\
            \Acu u(t,\om)
        :=
        \Big\{ \varphi \in C^{1,2}(\Lambda^t):
               \exists\;\H \in \Tc^t_+,~
               (\varphi- u^{t,\omega})_t(\0)
               =
               \displaystyle\max_{\tau \in \Tc^t}
               \E^{\P_{t,\om}}
               \big[ (\varphi - u^{t,\om})_{\tau \land \H} \big]
               \Big\}.
        \end{array}
        $$
The next definition requires the following notation for the
path-dependent second order differential operator on the shifted
canonical space: for all $(s,\om')\in\Lambda^t$,
        \b*
            (\Lc^{t,\om} \varphi )(s,\om')
            &:=&
            \partial_t \varphi(s, \om')
            + (\mu^{t,\om} \cdot \partial_{\om} \varphi)(s,\om')
            + \frac{1}{2} \big( (\sigma \sigma^T )^{t,\om} : \partial_{\om \om}^2 \varphi \big)(s,\om').
        \e*

    \begin{Definition} \label{def:visc_sol_PPDE}
    Let $u:\Lambda^0\longrightarrow\R$ be a locally bounded $\F^0-$progressive process.
    \\
        \rmi We say that $u \in \Ucb$ (resp. $u \in \Ucu$) is a viscosity subsolution (resp. supersolution) of PPDE \eqref{eq:PPDE} if, for any $(t,\om) \in [0,T) \x \Om^0$ and any $\varphi \in \Acb u(t, \om)$ (resp. $\varphi \in \Acu u(t, \om)$), it holds that
        \b*
            \big\{ - \Lc^{t,\om} \varphi
                   - \widehat{F}^{t,\om}(\cdot, u^{t,\om})
            \big\}(t,\0)
            &\le&
            (\mbox{resp.} \ge)~~ 0.
        \e*
        \rmii We say that $u$ is a viscosity solution of PPDE \eqref{eq:PPDE} if it is both a viscosity subsolution and a viscosity supersolution.
    \end{Definition}

    \begin{Remark} \label{rem:visc_sol_PPDE}
            \rmi In Definition \ref{def:visc_sol_PPDE}, we restrict ourselves, without loss of generality, to the test functions $\varphi \in \Acb$ (resp. $\Acu$) such that  $(\varphi- u^{t,\omega})_t(\0) = 0$.\\
        \rmii Similar to Remark 3.9 of Ekren, Keller, Touzi and Zhang \cite{EkrenKellerTouziZhang}, we can easily verify that under Assumption \ref{assum:fg}, for every $\lambda \in \R$, $u$ is a viscosity solution to \eqref{eq:PPDE} if and only if $\tilde u (t, \om) := e^{\lambda (T-t) } u(t,\om)$ is a viscosity solution of
        \b*
        - \Lc\tilde u
        -\widehat{F}_\lambda(., \tilde u)=0,
        &\mbox{where}&
        \widehat{F}_\lambda(t,\om,y)
        :=
        -\lambda y
        +e^{\lambda t}\widehat{F}\big(t,\om,e^{-\lambda t}y\big).
        \e*
        \rmiii Similar to Remark 2.11 of \cite{EkrenKellerTouziZhang}, we point out that in the Markovian setting, where the PPDE \eqref{eq:PPDE} reduces to a classical PDE, a viscosity solution in sense of Definition \ref{def:visc_sol_PPDE} is consistent to the viscosity solution in standard sense, by the uniqueness result proved below.
    \end{Remark}

\subsubsection{The existence and uniqueness of solutions to PPDE}

		This section follows closely the arguments of \cite{EkrenKellerTouziZhang, EkrenTouziZhang1, EkrenTouziZhang2}. However, their results do not apply to our context, because of the possible unboundedness of $\mu$ and $\sigma$. Moreover, the PPDE in our context is linear in the gradient and the Hessian components, which significantly simplifies the approach, see Remark 3.9 of \cite{EkrenTouziZhang1}.
		
    The above viscosity solution to PPDE \eqref{eq:PPDE} is closely related to the following decoupled FBSDE.
    For every $(t,\xb) \in \Lambda^0$, let $X^{0,t,\xb}$ be the solution of \eqref{eq:SDEt}, $(\widehat{Y}^{0,t,\xb}, \widehat{Z}^{0,t,\xb})$ be the solution of the BSDE on $(\Om^t, \Fc^t_T, \P_0^t)$,
    \be \label{eq:BSDE}
        \widehat{Y}_s
        &=&
        \xi(X_{\cdot}^{0, t,\xb})
        +
        \int_s^T \widehat{F} \big( r, X_{\cdot}^{0, t,\xb}, \widehat{Y}_r \big) dr
        - \int_s^T \widehat{Z}_r \cdot dB_r^t.
    \ee
By the Blumenthal 0-1 law, $\widehat{Y}^{0,t,\xb}_t$ is a constant and
we then define
    \be \label{eq:u0}
        \hat u(t,\xb) &:=& \widehat{Y}_t^{0, t,\xb}.
    \ee
For later use, we observe that, since the diffusion matrix $\sigma$
is nondegenerate, the above BSDE \eqref{eq:BSDE} is equivalent to
the following BSDE on $(\Om^t,\Fc^t_T,\P_{t,\xb})$:
    \b*
    \tilde Y_s
    &=&
    \xi^{t,\xb}(B^t_{\cdot})
    +
    \int_s^T \widehat{F}^{t,\xb} \big( r,B^t_{\cdot},
                                       \tilde Y_r\big)dr
                                      - \int_s^T \tilde Z_r \cdot
             \big(dB_r^t
                  -\mu^{t,\xb}(r,B^t_{\cdot})dr
             \big),
    \e*
    where $\widehat{F}^{t,\xb}$ is the shifted function of $\widehat{F}$ as introduced in \eqref{eq:shifted_proc}.
    Denote its solution by $(\tilde Y^{0, t,\xb}, \tilde Z^{0,t,\xb})$,
    then $\widehat{Y}_t^{0, t,\xb} = \tilde Y_t^{0, t,\xb} = \hat u(t,\xb)$ for every $(t,\xb)
\in \Lambda^0$. Moreover, by equation (4.6) of
\cite{EkrenTouziZhang1}, we have the dynamic programming principle
        \be \label{eq:DPP_BSDE}
            \tilde Y^{0,t,\xb}_s
            &=&
            \hat{u}^{t,\xb}(\tau, B^t_{\cdot})
            ~+~ \int_s^{\tau} \widehat{F}^{t,\xb}
              \big(r, B^t_{\cdot},\tilde Y_r^{0,t,\xb}\big)dr \nonumber \\
            &&-~
            \int_s^{\tau} \!\!\tilde Z_r^{0,t,\xb} \cdot
              \big( dB_r^t - \mu^{t,\xb}(r, B^t_{\cdot})dr
              \big),
                ~~\P_{t,\xb}-\mbox{a.s.},
        \ee
for all $(t,\xb) \in \Lambda^0$ and $\tau \in \Tc^t$.

    \vspace{2mm}

    Now, let us provide a representation of PPDE \eqref{eq:PPDE} by BSDE and a uniqueness result, whose proofs are very close to that in \cite{EkrenKellerTouziZhang, EkrenTouziZhang1, EkrenTouziZhang2},
    and we hence complete them in Appendix.

    \begin{Theorem} \label{theo:BSDE_PPDE}
        Let Assumption \ref{assum:fg} hold true.\\
        \rmi There is a constant $C>0$ such that $\forall (t,\om),(t',\om')\in\Lambda^0$,
        \b*
        |\hat{u}(t,\om)-\hat{u}(t',\om')|
        &\le&
        C \big(\| \om \|_t + \| \om' \|_{t'} \big) \big(\sqrt{|t-t'|}+\|\om_{t \land \cdot}-\om'_{t'\land \cdot}\|_T\big).
        \e*
        \rmii $\hat{u}$  is a viscosity solution to PPDE \eqref{eq:PPDE}.
    \end{Theorem}

	\begin{Theorem}\label{thm:comparison}
			Let Assumption \ref{assum:fg} hold true,
			$u^1, u^2$ be two $\F^0-$progressive c\`adl\`ag processes on $\Omega^0$ with corresponding jumps $\Delta u^1\ge 0\ge\Delta u^2$.
			Assume that $u^1$ (resp. $u^2$) is a viscosity subsolution (resp. supersolution) of PPDE \eqref{eq:PPDE},
			and $u^1(T,\cdot) \le \xi(\cdot) \le u^2(T,\cdot)$.
			Then $u^1 \le u^2$ on $\Lambda^0$.
    \end{Theorem}

\subsection{Proof of Theorem \ref{theo:Branching_BSDE}}
\label{subsec:proof_mainthm}

    Finally, we can complete the proof of our main result which gives a representation of BSDE by branching process.

    \vspace{2mm}

    \noindent {\bf Proof of Theorem \ref{theo:Branching_BSDE}}.
    By Theorems \ref{theo:BSDE_PPDE} and \ref{thm:comparison}, we only need to show that $v$ is a
    viscosity solution of \eqref{eq:PPDE} with terminal condition $\psi$ and generator
    $F_{n_0}$ defined in \eqref{eq:F_poly} following Definition \ref{def:visc_sol_PPDE}.
    We shall only show the subsolution part.
    Moreover, we recall that in the branching process,
    the process $^{t,\xb}X^{(1)}$ associated with the first particle is extended after its default time
    $T_1$ by $^{t,\xb}X^{(1)}_s := \Phi^{t,\xb}(s, B^{t,(1)})$ for all $s \in [t, T]$, where $B^{t,(1)}$ is defined by $B^{t,(1)}_{t+s} := W^1_s$ for all $s \in [t, T]$.

    Suppose that $v$ is not a viscosity subsolution of \eqref{eq:PPDE}, then by Definition \ref{def:visc_sol_PPDE} and Remark \ref{rem:visc_sol_PPDE}, there is
    $(t_0,\om_0) \in \Lambda^0$ and $\varphi \in \Acb v(t_0,\om_0)$ such that $v(t_0,\om_0)=\varphi(t_0,\0)$ and
    \be \label{eq:assum_contdic}
            - ~\Lcb \varphi (t_0, \om_0) ~=~ - ~\Lc \varphi (t_0, \om_0)
            ~-~ \beta \Gc \varphi (t_0, \om_0) &=& c ~>~ 0,
    \ee
    where $\Lc$ is defined by \eqref{eq:Lc} and
    \b*
            \Gc \varphi(t, \om) &:=&  \sum_{k=0}^{n_0} a_k(t,\om) \varphi^k(t,\om) - \varphi(t, \om) .
    \e*
    Without loss of generality, we suppose that $t_0 = 0$. Then $~^{0,\om_0}X^{(1)} = ~^{0,\0}X^{(1)}$.
    Further, it follows by the continuity of functions $\varphi$ and $v$ in Proposition \ref{prop:property_v} that
    for every $\eps > 0$,
    there is $\H \in \Tc_+^{0}$ such that for every $t \in [0, \bar \tau]$
    (with $\bar \tau := \H(~^{0,\0}X^{(1)}_{\cdot})$),
    \b*
            \big | v \big(t, ~^{0,\0}X^{(1)}_{\cdot} \big) - v \big(0, \0 \big) \big|
            +
            \big| \Gc \varphi \big(t,~^{0,\0}X^{(1)}_{\cdot}  \big)
            - \Gc v \big(t, ~^{0,\0}X^{(1)}_{\cdot} \big) e^{-\beta t} \big|
            &\le& \eps,
    \e*
    and
    \b*
        - ~\Lcb \varphi(t, ~^{0,\0}X^{(1)}_{\cdot})
        &\ge&
        c/2.
    \e*
    Clearly, $\bar \tau$ is a $\Fbb^{0}-$stopping time (see \eqref{eq:Fbb}) in probability space $(\Om, \Fc, \P)$.
    Denote $\H_h = \H \land h$ the $\F^0-$stopping time on $(\Om^0, \Fc^0, \P_0)$, $\taub_h := \bar \tau \land h$ and $X_t := ~^{0,\0}X^{(1)}_t$ the $\Fbb^0-$stopping time and process on $(\Om, \Fc, \P)$,
    it follows from equation \eqref{eq:DPP} of Lemma \ref{lemm:memoryless2}, together with \eqref{eq:assum_contdic},
    that
    \b*
						&&
                \E^{\P_X} \Big[ \varphi(\H_h, B^{0}_{\cdot}) - v(\H_h, B^{0}_{\cdot}) \Big] \\
            &=&
                \E^{\P} \Big[ \varphi(\taub_h, X_{\cdot}) - v(\taub_h, X_{\cdot}) \Big]\\
            &=&
                \E^{\P} \Big[ \varphi(\taub_h, X_{\cdot}) - \varphi(0, \0)
                + v(0, \0) -  v(\taub_h, X_{\cdot}) \Big] \\
            &=&
                \E^{\P} \Big[ \varphi(\taub_h, X_{\cdot}) - \varphi(0, \0) + \Gc v(T_1, X_{\cdot}) 1_{\taub_h \ge T_1} \Big]
                 + \E^{\P} \Big[ \big( v(T_1, X_{\cdot}) - v(\taub_h, X_{\cdot}) \big) 1_{\taub_h \ge T_1} \Big] \\
            &\le&
                 -~ \frac{c}{2} ~\E^{\P} \big[ \taub_h \big]
                 ~-~ \E \Big[ \int_{0}^{\taub_h} \Big( \beta \Gc \varphi(t, X_{\cdot})
                 ~-~ \Gc v(t,X_{\cdot}) \beta e^{- \beta t} \Big) dt ~\Big] \\
            &&
                  +~ \E^{\P} \Big[ \big( v(T_1, X_{\cdot}) - v(\taub_h, X_{\cdot}) \big) 1_{\taub_h \ge T_1} \Big] \\
            &\le&
                (-\frac{c}{2} + \eps \beta + 2 \eps) ~ \E^{\P}[ \taub_h ] ~~ < ~~ 0
    \e*
for $\eps$ small enough,
which is in contradiction with the fact that $\varphi \in \Acb v(t_0, \om_0)$ (see its definition below Assumption \ref{assum:fg}).
Therefore, $v$ is a viscosity subsolution of equation \eqref{eq:PPDE}. \qed

\section{Numerical examples}
\label{sec:num_example}

    In this section, we provide two numerical illustrations of our representation result, and the corresponding numerical implications.

\subsection{A two-dimensional example}

    Let us consider the following two decoupled FBSDEs:
    \be \label{eq:BSDE_example1}
        dX_t&=&\sigma X_t dB_t,~~~~~~~~~ X_0=1, \\
        dY_t&=&-\beta \big( F(Y_t)-Y_t \big)dt ~+~ Z_tdB_t,
    \ee
    with terminal condition $Y_T = \psi(X_T, A_T)$ and $A_t := \int_0^t X_s ds$, and the non-linearity $F$ is given by $F_1(y)=y^2$ or $F_2(y)=-y^2$.
    It is clear that the solution $Y$ can be given by the unique solution of PPDE
    \b*
        -~ \partial_t u (t,\om)
        ~-~ \frac{1}{2} \sigma^2 \om_t^2 \partial^2_{\om \om} u(t,\om)
        ~-~ F(u(t,\om))
        &=& 0,
    \e*
    with terminal condition $u(T,\om) := \psi(\om_T, \int_0^T \om_s ds)$.

    On the other hand, by adding a variable $a$, one can characterize the solution of FBSDE \eqref{eq:BSDE_example1} by some function $v(t,x,a)$ which is a classical viscosity solution of the following two PDEs:
\begin{eqnarray}
\partial_t v_1
+ x \partial_a v_1 + \frac{1}{2}\sigma^2 x^2 \partial^2_{xx} v_1 +
\beta (v_1^2-v_1)=0,\;
v_1(T,x,a)=\psi(x,a)  \;: \; \mathrm{PDE1}\\
\partial_t v_2
+ x \partial_a v_2 + \frac{1}{2} \sigma^2 x^2 \partial^2_{xx} v_2 +
\beta (-v_2^2-v_2)=0,\; v_2(T,x,a)=\psi(x,a)\;: \; \mathrm{PDE2}
\end{eqnarray}
These two-dimensional PDEs can be solved by a finite-difference method, which provide a benchmark for the evaluation of the performance of our Monte Carlo algorithm.

In our numerical experiments, we have taken a diffusion coefficient
$\sigma=0.2$ and a Poisson intensity  $\beta=0.1$, and the maturity
$T=2$ or $T=5$ years. For $T=2$ years (resp. $5$ years), the
probability of default is around $0.18$ (resp. $0.39$). The terminal condition is
$\psi(x,a)=(\frac{a}{T}-1)^+$.

In comparison with the KPP type PDE with $F_1(y)=y^2 $,
the replacement of the non-linearity $y^2$ by $-y^2$ has added the
term $(-1)^{N_T-1}$ in the multiplicative functional (see Equation
(\ref{eq:def_v})), without changing the complexity of the branching
diffusion  algorithm. More precisely, we have:
\begin{eqnarray}
{v}_1(0,X_0,A_0) &=& \E_{0,x}\Big[\prod_{i=1}^{N_T}
\psi(X_T^i,A_T^i)\Big] ,
\\
{v}_2(0,X_0,A_0) &=& \E_{0,x}\Big[ (-1)^{N_T-1} \prod_{i=1}^{N_T}
\psi(X_T^i,A_T^i)
        \Big] .
\end{eqnarray}

Our branching diffusion algorithm has been checked against
a two-dimensional PDE solver with an ADI scheme (see Tables
\ref{BMMTest1}, \ref{BMMTest2}). The degenerate PDEs have been
converted into elliptic PDEs by introducing the process
$\tilde{A}_t=\int_0^t X_s ds +(T-t)X_t$, satisfying
$d\tilde{A}_t=(T-t)dX_t$. The computational experiments was done
using a PC with 2.99 Ghz Intel Core 2 Duo CPU.
Note that our algorithm converges to the exact PDE result
as expected and the error is properly indicated by the Monte-Carlo
standard deviation estimator (see column Stdev). In order to
illustrate the impact of the non-linearity $F$ on the price $v$, we
have indicated the price corresponding to $\beta=0$.

\begin{table}[h]\begin{center}
\begin{tabular}{||c|c|c|c|c|c||}
\hline\hline N & Fair(PDE1) & Stdev(PDE1) & Fair(PDE2) & Stdev(PDE2)
&  CPU (seconds)
\\ \hline

$12$ & $5.69$ & $0.16$ & $5.36$ & $0.16$ & $0.1$  \\ \hline

$14$ & $5.61$ & $0.08$ & $5.23$ & $0.08$ & $0.6$ \\ \hline

$16$ & $5.50$ & $0.04$ & $5.15$ & $0.04$ & $1.5$ \\ \hline

$18$ & $5.52$ & $0.02$ & $5.16$ & $0.02$ & $5.9$ \\ \hline

$20$ & $5.53$ & $0.01$ & $5.16$ & $0.01$ & $23.6$ \\ \hline

$22$ & $5.54$ & $0.00$ & $5.17$ & $0.01$ & $94.1$ \\ \hline\hline
\end{tabular}
\end{center}
\caption{MC price quoted in percent as a function of the number of
MC paths $2^N$. PDE pricer(PDE1) = ${\bf 5.54}$. PDE pricer(PDE2) =
${\bf 5.17}$ (CPU PDE: {\bf $10$ seconds}). Maturity= $2$ years.
Non-linearities for PDE1 (resp. PDE2) $F_1(u)=u^2$ (resp.
$F_2(u)=-u^2$). For completeness, the price with $\beta=0$ (which
can be obtained using a classical Monte-Carlo pricer) is $6.52$.}
\label{BMMTest1}
\end{table}

\begin{table}[h]\begin{center}
\begin{tabular}{||c|c|c|c|c|c||}
\hline\hline N & Fair(PDE1) & Stdev(PDE1) & Fair(PDE2) & Stdev(PDE2)
& CPU (seconds)
\\ \hline

$12$ & $7.40$ & $0.25$ & $5.63$ & $0.26$ & $0.3$\\ \hline

$14$ & $7.28$ & $0.12$ & $5.60$ & $0.13$ & $1.1$\\ \hline

$16$ & $7.20$ & $0.06$ & $5.47$ & $0.07$ & $4.3$\\ \hline

$18$ & $7.24$ & $0.03$ & $5.48$ & $0.03$ & $17.0$ \\ \hline

$20$ & $7.24$ & $0.02$ & $5.50$ & $0.02$ & $68.3$ \\\hline

$22$ & $7.24$ & $0.01$ & $5.51$ & $0.01$ & $272.9$ \\ \hline\hline

\end{tabular}
\end{center}
\caption{MC price quoted in percent as a function of the number of
MC paths $2^N$. PDE pricer(PDE1) = ${\bf 7.24}$. PDE pricer(PDE2) =
${\bf 5.51}$ (CPU PDE: {\bf $25$ seconds}). Maturity= $5$ years.
Non-linearities for PDE1 (resp. PDE2) $F_1(u)=u^2$ (resp.
$F_2(u)=-u^2$). For completeness, the price with $\beta=0$ (which
can be obtained using a classical Monte-Carlo pricer) is $10.24$.}
\label{BMMTest2}
\end{table}

\subsection{An eight-dimensional example}

We would like to highlight that the high-dimensional case can be easily handled in our framework 
by simulating the branching particles with a high-dimensional diffusion process. 
This is out-of-reach with finite-difference scheme methods and not such an easy step for the classical numerical schemes of BSDEs which require computing conditional expectations. 
In order to illustrate this point, we have implemented our algorithm for the following  decoupled FBSDEs
    \be
        dX_t^i&=&\sigma_i X_t^i dB_t^i,\quad  d\langle B^i,B^j\rangle_t=\delta_{i,j}dt, \quad X_0^i=1,\quad i,j=1,\ldots,4, \\
        dY_t&=&-\beta \big( F(Y_t)-Y_t \big)dt ~+~ \sum_{i=1}^4 Z_t^idB_t^i,
    \ee
    with terminal condition $Y_T = \psi(X_T, A_T)$, $A_t^i := \int_0^t X_s^i ds$, and the non-linearity $F$ is given by $F_1(y)=y^2$ or $F_2(y)=-y^2$.  $X_t=(X_t^i)_{i=1,\ldots ,4}$ define a $4$d uncorrelated
    geometric Brownian motion and we have $4$ path-dependent variables $A_t=(A_t^i)_{i=1,\ldots ,4}$. 
		Similarly, the solution is related to the non-linear $8$d-PDEs
\begin{eqnarray}
\partial_t v_1
+ {\cal L} v_1 +
\beta (v_1^2-v_1)=0,\;
v_1(T,x,a)=\psi(x,a)  \;: \; \mathrm{PDE1}\\
\partial_t v_2
+{\cal L} v_2 +
\beta (-v_2^2-v_2)=0,\; v_2(T,x,a)=\psi(x,a)\;: \; \mathrm{PDE2}
\end{eqnarray}
with  ${\cal L}=\frac{1}{2}\sum_{i=1}^4 \sigma_i^2 \partial_{x_i}^2 +
\sum_{i=1}^4 x_i \partial_{a_i}$. In our numerical experiments, we have taken a diffusion coefficient
$\sigma_i=0.2$, a Poisson intensity  $\beta=0.1$, and the maturity
$T=2$ or $T=5$ years. The terminal condition is  $\psi(x,a)=(\frac{\sum_{i=1}^4 a_i}{4T}-1)^+$.

These eight-dimensional PDEs  suffer from the curse of dimensionality and we are unable to
 solve them by a finite-difference method. In the particular case of
  a constant terminal condition $v_1(T,x,a)=v_2(T,x,a)=1/2$, these PDEs reduce to ODEs which can
  be integrated out explicitly: $v_1(0,X_0,0)^{-1}=1+e^{\beta T}$,
$v_2(0,X_0,0)^{-1}=-1+3e^{\beta T}$. As a simple preliminary benchmark, we have checked that our numerical algorithm
 reproduces exactly these solutions.  In the case of the non-trivial payoff $\psi(x,a)=(\frac{\sum_{i=1}^4 a_i}{4T}-1)^+$,
  we have checked that our branching
diffusion algorithm converge  (see Tables
\ref{BMMTest3}, \ref{BMMTest4}). 
We also report the average number of descendants
generated up to the maturity $T$.  As far as we know, we are not unaware of
alternative numerical methods for solving such a non-linear 8d-PDE.

\begin{table}[h]\begin{center}
\begin{tabular}{||c|c|c|c|c|c||}
\hline\hline N & Fair(PDE1) & Stdev(PDE1) & Fair(PDE2) & Stdev(PDE2)
\\ \hline

$12$ & $2.77$ & $0.08$ & $2.67$ & $0.08$   \\ \hline

$14$ & $2.69$ & $0.04$ & $2.60$ & $0.04$   \\ \hline

$16$ & $2.71$ & $0.02$ & $2.62$ & $0.02$   \\ \hline

$18$ & $2.72$ & $0.01$ & $2.63$ & $0.01$   \\ \hline

$20$ & $2.74$ & $0.00$ & $2.65$ & $0.00$   \\ \hline

$22$ & $2.74$ & $0.00$ & $2.65$ & $0.00$   \\ \hline

\end{tabular}
\end{center}
\caption{MC price quoted in percent as a function of the number of
MC paths $2^N$. Maturity= $2$ years.
Non-linearities for PDE1 (resp. PDE2) $F_1(u)=u^2$ (resp.
$F_2(u)=-u^2$). For completeness, the price with $\beta=0$ (which
can be obtained using a classical Monte-Carlo pricer) is $3.29$. The average number of descendants generated is
$1.22$. }
\label{BMMTest3}
\end{table}

\begin{table}[h]\begin{center}
\begin{tabular}{||c|c|c|c|c|c||}
\hline\hline N & Fair(PDE1) & Stdev(PDE1) & Fair(PDE2) & Stdev(PDE2)

\\ \hline

$12$ & $3.35$ & $0.11$ & $2.99$ & $0.11$  \\ \hline

$14$ & $3.40$ & $0.06$ & $3.04$ & $0.06$  \\ \hline

$16$ & $3.38$ & $0.03$ & $3.01$ & $0.03$  \\ \hline

$18$ & $3.38$ & $0.01$ & $2.99$ & $0.01$   \\ \hline

$20$ & $3.38$ & $0.01$ &  $3.00$ & $0.01$     \\ \hline

$22$ & $3.38$ & $0.00$ &  $3.00$ & $0.00$    \\ \hline

\hline

\end{tabular}
\end{center}
\caption{MC price quoted in percent as a function of the number of
MC paths $2^N$. Maturity= $5$ years.
Non-linearities for PDE1 (resp. PDE2) $F_1(u)=u^2$ (resp.
$F_2(u)=-u^2$). For completeness, the price with $\beta=0$ (which
can be obtained using a classical Monte-Carlo pricer) is $5.24$. The average number of descendants generated is $1.65$.}
\label{BMMTest4}
\end{table}

\appendix

\section{Appendix}

    Here we complete the proofs for Theorems \ref{theo:BSDE_PPDE} and \ref{thm:comparison}, where the arguments are mainly adapted from that in Ekren, Touzi and Zhang \cite{EkrenTouziZhang1, EkrenTouziZhang2}.

    \vspace{2mm}

    \noindent {\bf Proof of Theorem \ref{theo:BSDE_PPDE}.}
    \rmi is proved in Proposition 4.5 of \cite{EkrenTouziZhang1}, since our BSDE \eqref{eq:BSDE} is a
particular case to their equation (4.4).
        It is in fact an immediate consequence of Proposition 2.1 in El Karoui, Peng and Quenez \cite{ElKarouiBSDE}
         together with the estimation in our Lemma \ref{lemm:SDE_Lip}.
\\
        \rmii We adapt the arguments in the proof of Proposition 4.5 \cite{EkrenTouziZhang1} to our context.
        We only show that $\hat{u}$ is a viscosity subsolution. Assume $\hat{u}$ is not a viscosity subsolution,
        then there exist $(t,\om) \in \Lambda^0$ and $\varphi \in \Acb \hat{u}(t,\om)$ such that
        $$
            c
            ~:=~
            - \Lc^{t,\om}\varphi(t,\0) - \widehat{F}^{t,\om}(t, \0, \hat{u}(t, \om))
            ~>~
             0.
        $$
		Without loss of generality, we may also assume that
        $\varphi(t,\0)=\hat u(t,\om)$.
        Denote, for $s \in [t,T]$,
        \b*
                \tilde Y'_s := \varphi(s, B^t),
                ~~ \tilde Z'_s := \partial_{\om} \varphi(s, B^t),
                ~~ \delta Y_s := \tilde Y'_s - \tilde Y_s^{0,t,\om},
                ~~ \delta Z_s := \tilde Z'_s - \tilde Z_s^{0,t,\om}.
        \e*
        Applying It\^o's formula, we have $\P_{t,\om}-a.s$,
        \b*
                d(\delta Y_s)
                &=&
                \Big[ (\Lc^{t,\om} \varphi)(s, B^t_{\cdot})
                    +
                    \widehat{F}^{t,\om} \big(s, B^t_{\cdot}, \tilde Y_s^{0,t,\om} \big) \Big] ds
                + \delta Z_s \cdot \big( d B_s^t - \mu^{t,\om}(s, B^t_{\cdot})  ds \big) \\
                &=&
                \Big[ (\Lc^{t,\om} \varphi)(s, B^t_{\cdot})
                    +
                    \widehat{F}^{t,\om} \big(s, B^t_{\cdot}, \tilde Y_s' \big)
                    +\alpha_s \delta Y_s \Big] ds
                + \delta Z_s \cdot \big( d B_s^t - \mu^{t,\om}(s, B^t_{\cdot})  ds \big),
        \e*
        where $\alpha$ is a $\F^t$-progressively measurable process bounded by the Lipschitz constant $L_0$ of $\widehat F$ in $y$.

        Observing that $\tilde Y_t' = \varphi(t,\0) = \hat u(t,\om)$ and $\delta Y_t = 0$, we define a stopping time
        \b*
                \H
                &:=&
                T \land
                \inf \big\{ s > t ~:
                        - \Lc^{t,\om} \varphi(s, B^t_{\cdot}) - \widehat{F}^{t,\om}(s, B^t_{\cdot}, \varphi(s, B^t_{\cdot}))
                        - L_0 |\delta Y_s |
                        ~\le~
                        c /2
                     \big\}.
        \e*
        Then by the continuity of $\Lc^{t,\om} \varphi$ as well as $ \widehat{F} $, we have $\H \in \Tc^t_+$ and
        \b*
            - \Lc^{t,\om} \varphi (s, B^t_{\cdot})
            - \widehat{F}^{t,\om}(s, B^t_{\cdot}, \tilde Y_s')
            - \alpha_s \delta Y_s
            ~\ge~ c/2, && \mbox{for all}~ s \in [t, \H].
        \e*
        Now for any $\tau \in \Tc^t$ such that $\tau \le \H$, we have
        \be \label{eq:deltaY}
                0 &=& \delta Y_t
                        ~=~
                            \delta Y_{\tau} - \int_t^{\tau} \Big[ (\Lc^{t,\om}\varphi)(s, B^t_{\cdot})
                            + \widehat{F}^{t,\om}(s, B^t_{\cdot}, \tilde Y'_s)
                            + \alpha_s \delta Y_s \Big] \nonumber \\
                &&~~~~~~~~~~~ - \int_t^{\tau} \delta Z_s \cdot
                        \big( d B_s^t - \mu^{t,\om}(s, B^t_{\cdot})  ds \big) \nonumber \\
                &\ge&
                            \varphi(\tau, B^t) - \hat u^{t,\om}(\tau, B^t)
                            + c(\tau-t)/2
                            - \int_t^{\tau} \delta Z_s \cdot \big( d B_s^t - \mu^{t,\om}(s, B^t_{\cdot})  ds \big),
        \ee
        $\P_{t,\om}-a.s.$.
        We recall that $\P_{t, \om}$ is defined by \eqref{eq:Pt}, under which the canonical process $B^t$ is a solution to SDE \eqref{eq:SDEt}.
        Therefore,
        \b*
            \int_t^{\cdot} dB^t_s - \mu^{t,\om}(s, B^t_{\cdot}) ds
        \e*
        is a $\P_{t,\om}$-martingale.
        By taking expectation on \eqref{eq:deltaY} under $\P_{t,\om}$, it follows that $\E^{\P_{t,\om}}_t[ (\varphi -\hat{u}^{t,\om})(\tau, B^t_{\cdot}) ] < 0$,
        which contradicts the fact that $\varphi \in \Acb\hat{u}(t,\om)$.
        \qed

\vspace{5mm}

    In preparation of the comparison principle in Theorem \ref{thm:comparison}, we first introduce two extended spaces $\Cu^{1,2}_{t,\om}(\Lambda^t)$ and $\Cb^{1,2}_{t,\om}(\Lambda^t)$ of $C^{1,2}(\Lambda^0)$ and derive a partial comparison principle as in \cite{EkrenTouziZhang1, EkrenTouziZhang2}.

    \begin{Definition}\label{defC12bar}
        Let $(t,\om) \in \Lambda^0$, $u: \Lambda^t \to \R$ be $\F^t-$adapted.\\
        \rmi We say $u \in \Cu^{1,2}_{t,\om}(\Lambda^t)$ if there exist an increasing sequence of $\F^t-$stopping times in $\Tc^t$: $t = \tau_0 \le \tau_1 \le \cdots \le T$ such that,\\
        \rma $\tau_i < \tau_{i+1}$ whenever $\tau_i < T$, and for all $\omt \in \Om^t$, the set $\{i: \tau_i(\omt) < T \}$ is finite; \\
        \rmb For each $i \ge 0$ and $\omt \in \Om^t$, $\tau^{\tau_i(\omt), \omt}_{i+1} \in \Tc^{\tau_i(\omt)}$
        and $u^{\tau_i(\omt),\omt} \in C_b^{1,2} \big( \Lambda^{\tau_i(\omt)} \big(\tau_{i+1}^{\tau_i(\omt),\omt} \big) \big)$;
        \rmc $u$ has non-negative jumps ($\Delta u\ge 0$), and
        \be \label{eq:cond_Cb12}
        \E^{\P_{t,\om}}
        \Big[ \sum_{i\ge 0} \int_{\tau_i}^{\tau_{i+1}}
              \Big( \big|\Lc^{t,\om} u \big|^2
                    +\big|\sigma^{t,\om}\partial_{\om}u\big|^2
              \Big) (s, B^t) ds \Big]
        <\infty.
        \ee
        \rmii We say $u \in \Cb^{1,2}_{t,\om}(\Lambda^t)$ if $-u \in \Cu^{1,2}_{t,\om}(\Lambda^t)$.
    \end{Definition}

        \begin{Lemma}\label{lemm:partial_comp}
            Suppose that Assumption \ref{assum:fg} holds true. Let $u^1$ be a viscosity subsolution and $u^2$ be a viscosity supersolution of PPDE \eqref{eq:PPDE} such that $u^1(T, \cdot) \le u^2(T,\cdot)$.
            If  $u^1 \in \Cu^{1,2}_{0,0}(\Lambda^0)$ or $u^2 \in \Cb^{1,2}_{0,0}(\Lambda^0)$, then $u^1 \le u^2$ on $\Lambda^0$.
        \end{Lemma}
        \proof We follow the lines of Proposition 4.1 of Ekren, Touzi and Zhang \cite{EkrenTouziZhang2}.
        Suppose that $u^1 \in \Cu^{1,2}_{0,0}(\Lambda^0)$.
        First, let us show that, for every $i \ge 0$ and $\om \in \Om^0$,
        \be \label{eq:PC_iteration}
                \big( u^1 - u^2 \big)^+_{\tau_i(\om)}(\om)
                &\le&
                \E^{\P_{\tau_i(\om),\om}} \Big[ \Big( (u^1)^{\tau_i(\om),\om}_{\tau_{i+1}(\om)} - (u^2)^{\tau_i(\om),\om}_{\tau_{i+1}(\om)} \Big)^+ \Big].
        \ee
        Without loss of generality, it is enough to consider the case $i = 0$, where $\P_{\tau_0(\om), \om} = \P_X$ for all $\om \in \Om$. Assume to the contrary that
        \b*
            2Tc
            &:=& (u^1_0 - u^2_0)^+
            ~-~ \E^{\P_X} \big[ \big( u^1_{\tau_1} - u^2_{\tau_1} \big)^+ \big]
            ~>~ 0,
        \e*
        we set
        \b*
            X_t := (u^1_t - u^2_t)^+ +  c t,
            & Y_t ~:= \sup_{\tau \in \Tc^t} \E_t [X_{\tau \land \tau_1}],
            & \tau^* := \inf \{ t > 0 : X_t = Y_t \} \le \tau_1,
        \e*
        where $\E_t[ \zeta](\om) := \E^{\P_{t,\om}}[\zeta^{t,\om}] = \E^{\P_X}[ \zeta | \Fc_t](\om)$. We notice that the conditional expectation $\E_t$ is defined by using shifting operators, and in this case the supremum in the definition of $Y$ is the same of the essential supremum (see also Theorem 2.3 of Nutz and van Handel \cite{NutzHandel} for details for a similar problem).
        In particular, $\E_0[\cdot] = \E^{\P_X}[\cdot]$.
        Then $(Y_t)_{t \ge 0}$ is a supermartingale, $(Y_{t \land \tau^*})_{t \ge 0 }$ is a martingale and $\tau^*$ is an optimal stopping time for the problem $\sup_{\tau \in \Tc^0} \E_0 [ X_{\tau}]$. It follows that
        \b*
                \E_0 [ X_{\tau^*}]
                ~=~ \E_0 [ Y_{\tau^*}]
                ~=~ Y_0
                ~\ge~ X_0
                &=~ 2Tc + \E_t \big[ \big( u^1_{\tau_1} - u^2_{\tau_1} \big)^+ \big]
                ~\ge~ Tc + \E_0 [ X_{\tau_1}].
        \e*
        Then there exists $\om^* \in \Om^0$ such that $t^* := \tau^*(\om^*) < \tau_1$. And therefore
        \b*
                (u^1 - u^2)^+(t^*, \om^*) + c t^* = X_{t^*}(\om^*) = Y_{t^*}(\om^*) \ge \E_{t^*}\big[ X_{\tau_1} \big] > c t^*,
        \e*
        which implies that $0 < (u^1 - u^2)^+(t^*, \om^*)$.
        Set $\varphi(t,\om) := (u^1)^{t^*, \om^*}(t,\om) + c(t^*)$.
        Then $\varphi \in C^{1,2}(\Lambda^{t^*}(\tau_1))$ since $u^1 \in C^{1,2}(\Lambda(\tau_1))$. Moreover, let
        \b*
                \H &:=& \inf  \big\{ t > t^* ~:~ u^1_t - u^2_t \le 0 \big \} \land \tau_1 ~\in~ \Tc^{t^*}_+ .
        \e*
        Then for every $\tau \in \Tc^{t^*}$,
        \b*
                (\varphi - (u^2)^{t^*, \om^*})(t^*,\0)
                &=& X_{t^*}(\om^*)
                ~\ge~
                \E_{t^*} \big[ Y_{\tau \land \H} \big](\om^*) \\
                &\ge& \E_{t^*} \big[ X_{\tau \land \H} \big](\om^*)
                ~=~ \E^{\P_{t^*,\om^*}} \big[ \big( \varphi - (u^2)^{t^*, \om^*} \big)_{\tau \land \H} \big],
        \e*
        which implies that $\varphi \in \Acb u^2(t^*, \om^*)$. It follows that
        \b*
                0
                &\le&
                \big\{ ~-~ \Lc \varphi - \widehat F(\cdot, \varphi) \big\}(t^*, \om^*)
                ~\le~
                -~ c
                - \big(~-~ \Lc u^1 - \widehat F(\cdot, u^1) \big) (t^*, \om^*) ,
        \e*
        which contradicts the fact that $u^1$ is a subsolution and we hence prove \eqref{eq:PC_iteration}.
        Further, since $(\P_{\tau_i(\om), \om})_{\om \in \Om}$ induces a r.c.p.d. of $\P_X$ w.r.t. $\Fc_{\tau_i}$ (see Remark \ref{rem:rcpd}), it follows by \eqref{eq:PC_iteration} that for every $ i \ge 0$,
        \b*
                (u^1 - u^2)_0
                &\le&
                \E_0 \big[ (u^1 - u^2)^+_{\tau_i} \big].
        \e*
        By sending $i \to \infty$, we get that $(u^1 - u^2)_0 \le \E_0 [ (u^1 - u^2)_T^+ ] = 0$,
        which completes the proof of $u_0^1 \le u_0^2$. \qed

        \vspace{2mm}

        \noindent {\bf Proof of Theorem \ref{thm:comparison}}. We follows the lines of the proof of Theorem 7.4 of Ekren, Touzi and Zhang \cite{EkrenTouziZhang1}, where a comparison principle for PPDE \eqref{eq:PPDE} was proved in case $\sigma \equiv I_d$. In spirit of Remark \ref{rem:visc_sol_PPDE}, we suppose without loss of generality that $\widehat F$ decreases in $y$.

        For every $\eps > 0$, we denote
        \b*
            O_{\eps} := \{ x \in \R^d ~: |x| < \eps \},
            & \Ob_{\eps} := \{ x \in \R^d ~: |x| \le \eps \},
            & \partial O_{\eps} := \{ x \in \R^d ~: |x| = \eps \};
        \e*
        \b*
            \Oc^{\eps}_t := [t,T) \x O_{\eps},
            ~~ \Ocb^{\eps}_t := [t,T] \x \Ob_{\eps},
            ~~ \partial \Oc^{\eps}_t := ([t,T] \x \partial O_{\eps}) \cup (\{T \} \x O_{\eps}).
        \e*
        Let $t_0 = 0$, $x_0 = 0$, $(t_i)_{i \ge 1}$ an increasing sequence in $(0, T]$ with $t_i = T$ when $i$ is large enough, and $(x_i)_{i \ge 1}$ a sequence in $\R^d$. Set $\pi := (t_i, x_i)_{i \ge 0}$ and $\pi_n := (t_i, x_i)_{0 \le i \le n}$. Given $\pi_n$ and $(t,x) \in \Oc^{\eps}_{t_n}$, define
        \b*
            \H_0^{t,x,\eps} := \inf \{ s \ge t : |B^t +x| = \eps \} \land T, & \H_{i+1}^{t,x,\eps} := \inf \{ s \ge \H_i^{t,x,\eps} : |B^t_s - B^t_{\H_i^{t,x,\eps}}| = \eps \} \land T.
        \e*
        For $t \in (t_n, T]$, let $\Bh^{\eps, \pi_n, t, x}(\om)$ denote the continuous path on $[0,T]$ obtained by linear interpolation of the function $b(t_i) := \sum_{j=0}^i x_j$ for $0 \le i \le n$ and
        $b(\H_i^{t,x,\eps}) := \sum_{j=0}^n x_j + x + B^t_{\H_i^{t,x,\eps}}(\om)$ for all $i \ge 0$.
        Define
        \b*
            \theta_n^{\eps}(\pi_n; (t,x) ) &:=& \Yc_t^{\eps, \pi_n, t,x},
        \e*
        where, omitting the superscripts $^{\eps, \pi_n, t,x}$, $\Yc$ is defined under $\P_{t,\Bh}$ by
        \b*
        \Yc_s
        &=&
        \xi(\Bh)
        +
        \int_s^T
        \widehat F \Big(r, \sum_{i \ge -1}
                   \Bh_{\cdot \land \H_i^{t,x,\eps}}
                   1_{[\H_i^{t,x,\eps}, \H_{i+1}^{t,x,\eps})}
                 ,\Yc_r
           \Big) dr
        - \int_s^T \Zc_r \big(dB_r-\mu(r,\Bh)dr\big),
        \e*
with $\H_{-1}^{t,x,\eps} := t$. Then clearly, for every $n$ and
$\pi_n$, the deterministic function $\theta_n^{\eps} :=
\theta_n^{\eps}(\pi_n; \cdot)$ is the viscosity solution of the
standard PDE on $\Oc^{\eps}_{t_n}$:
        \begin{equation} \label{eq:PDE_O}
            -\partial_t \theta^{\eps}_n
            - \mu(s, \omh^{\pi_n}) D \theta_n^{\eps}
            - \frac{1}{2} \sigma \sigma^T(s, \omh^{\pi_n}): D^2 \theta^{\eps}_n
            - \widehat F(s, \omh^{\pi_n}, \theta_n^{\eps})
            = 0
            ~\mbox{on}~ \Oc^{\eps}_{t_n},
        \end{equation}
        with terminal condition $\theta^{\eps}_n(\pi_n; t,x) = \theta^{\eps}_{n+1}(\pi_n, (t,x); t, 0)$ on $\partial \Oc^{\eps}_{t_n}$, where $\omh^{\pi_n} := \Bh^{\eps, \pi_n, t, x}_{\cdot \land t_n}$ is deterministic, and $\theta^{\eps}_n(\pi_n; T,x) = \xi(\omh^{\pi_n})$ when $t_n = T$.
        Further, since $\sigma \sigma^T$ is non-degenerate, it follows from Proposition 7.2 of \cite{EkrenTouziZhang1} that for every $\delta > 0$, there is $\thetau^{\eps, \delta}_n \in C^{1,2}(\Oc^{\eps}_{t_n})$ which is a classical supersolution of \eqref{eq:PDE_O} such that $\thetau^{\eps, \delta}_n(\pi_n; t,x) \ge \theta^{\eps}_{n+1}(\pi_n, (t,x); t, 0)$ on $\partial \Oc^{\eps}_{t_n}$ and $|\thetau^{\eps,\delta}_n - \theta^{\eps}_n| \le \delta $ on $\Ocb^{\eps}_{t_n}$. Let $\delta_n = \eps/2^n$, $\H^{\eps}_i := \H_i^{0,0,\eps}$, and $\Bh^{\eps}$ be the linear interpolation of $(\H_i^{\eps}, B_{\H_i^{\eps}})_{i \ge 0}$. Define
        \b*
            \psi^{\eps}(t,\om)
            &:=&
            \sum_{n=0}^{\infty} \Big( \delta_n + \thetau_n^{\eps, \delta_n} \big( (\H_i^{\eps}, B_{\H_i^{\eps}})_{0 \le i \le n}; t, B_t - B_{\H_n^{\eps}} \big) \Big) 1_{[\H_n^{\eps}, \H_{n+1}^{\eps})},
        \e*
        and denote
        \b*
                \Bt^{\eps}_{\cdot} &:=& \sum_{i \ge -1} \Bh^{\eps}_{\cdot \land \H_i^{t,x,\eps}} 1_{[H_i^{t,x,\eps}, H_{i+1}^{t,x,\eps})}.
        \e*
        One can check straightforwardly that $-\psi^{\eps}$ satisfies the conditions of Definition \ref{defC12bar} (c), $\psi^{\eps}(T,\om) \ge \xi(\Bt^{\eps})$, and
        \be
            - \partial_t \psi^{\eps} - \mu(s, \Bt^{\eps}_{\cdot} ) \cdot \partial_{\om} \psi^{\eps} - \sigma(s, \Bt^{\eps}_{\cdot} ) : \partial_{\om\om} \psi^{\eps} - \widehat F \Big(s, \Bt^{\eps}_{\cdot}, \psi^{\eps}(s, \Bt^{\eps}) \Big) \ge 0.
        \ee
        Then $\Yt := \psi^{\eps}$, $\Zt := \partial_{\om} \psi^{\eps}$ satisfy the BSDE
        \b*
            \Yt_s
            &=&
            \Yt_{\H_{i+1}}
            ~+~ \int_s^T
            \widehat F \Big(r, \Bt^{\eps}_{\cdot} , \Yt_r \Big) dr
            - \int_s^T \Zt_r \Big( dB_r - \mu \Big(r, \Bt_{\cdot}^{\eps} \Big) dr \Big), ~~ \P_X-a.s.
        \e*
        on every interval $[\H_i, \H_{i+1})$ such that $ \sup_{0 \le t \le T} |\Yc_t - \Yt_t| \le \eps$, which implies that \eqref{eq:cond_Cb12} holds true for $\psi^{\eps}$ and hence $\psi^{\eps} \in \Cb_{0,0}^{1,2}(\Lambda^0)$.
        Notice that $\|\Bt^{\eps} - B \|_T \le \eps$, then for some constant $C$,
        \b*
                | \xi(\Bt^{\eps}) - \xi(B) | ~\le~ C \eps,
                &&
                |\widehat F (s, \Bt^{\eps},y) - \widehat F(s, B, y)| ~\le~ C \eps.
        \e*
        Set
        \b*
            \psiu := \psi^{\eps} + 2 C \eps [1 + T-t],
        \e*
        one can verify that $\psiu \in \Cb^{1,2}_{0,0}(\Lambda^0)$ is a viscosity supersolution of \eqref{eq:PPDE},
        and it follows by the partial comparison principle in Lemma \ref{lemm:partial_comp} that $u^1(0,\0) \le \psiu(0,\0)$. Similarly, we can construct a viscosity subsolution $\psib \in \underline{C}^{1,2}_{0,\0}(\Lambda^0)$ such that
        $u^2(0, \0) \ge \psib(0,\0)$ and $|\psiu - \psib| \le 4 C[2+T] \eps$.
        By sending $\eps \to 0$, we conclude the proof.
        \qed

\end{document}